\documentclass[11pt]{article}
\usepackage[margin=1in]{geometry}  
\usepackage{graphicx}
\usepackage{subcaption}              
\usepackage{amsmath}               
\usepackage{amsfonts}              
\usepackage{amsthm}
\usepackage{cancel, soul, xcolor}
\usepackage[square, comma, numbers, sort & compress]{natbib} 
\usepackage{hyperref}
\usepackage{mathtools,breqn}
\usepackage{tikz}
\usetikzlibrary{decorations.pathmorphing}
\usetikzlibrary{shapes}
\usetikzlibrary{arrows,decorations.pathreplacing}
\usepackage{float}
 \usepackage[titletoc,toc,title]{appendix}

\allowdisplaybreaks
\usepackage{resizegather}

\usepackage{epstopdf}

\newtheorem{Def}{Definition}
\newtheorem{Thm}{Theorem}

\newtheorem{Prop}{Proposition}
\newtheorem*{Pf}{Proof}
\newtheorem{Lem}{Lemma}
\newtheorem{Rk}{Remark}
\newtheorem{Cc}{Corollary}

\setstcolor{red}

\begin{document}

\title{An analysis of the input-to-state-stabilisation of linear hyperbolic systems of balance laws with boundary disturbances}

\author{Gediyon Y. Weldegiyorgis, Mapundi K. Banda}

\maketitle

\begin{abstract}
In this paper, a linear hyperbolic system of balance laws with boundary disturbances in one dimension is considered. An explicit candidate Input-to-State Stability (ISS)-Lyapunov function in $ L^2- $norm is considered and discretised to investigate conditions for ISS of the discrete system as well. Finally, experimental results on test examples including the Saint-Venant equations with boundary disturbances are presented. The numerical results demonstrate the expected theoretical decay of the Lyapunov function.   
\end{abstract}





\section{Introduction}\label{sec:intro}
We consider a $ k\times k$ system described by the following linear hyperbolic system of balance laws with variable coefficients
\begin{equation}\label{eq:LHSBLaws}
\partial_t W(x,t) + \Lambda(x) \partial_x W(x,t) + \Pi(x) W(x,t) = 0, \quad (x,t) \in [0,l]\times [0, +\infty), 
\end{equation}
where $ W := W(x,t): [0,l]\times [0, +\infty) \rightarrow \mathbb{R}^k $ is a state vector, $\Lambda(x) = \text{diag}\{{\Lambda^+}(x), -{\Lambda^-}(x) \} $, with $ {\Lambda^+}(x) \in \mathbb{R}_{+}^{m \times m}$ and $ {\Lambda^-}(x) \in \mathbb{R}_{+}^{(k-m) \times (k-m)}$, is a non-zero diagonal matrix and $ \Pi(x) \in \mathbb{R}^{k\times k} $ is a non-zero matrix. By using the diagonal entries of $ \Lambda(x) $, the state vector $ W $ is specified by $ W = [{W^+}, {W^-}]^\top $, where $ {W^+} \in \mathbb{R}^{m}$ and $ {W^-} \in \mathbb{R}^{k-m}$. 

The system \eqref{eq:LHSBLaws} is subject to an initial condition set as
\begin{equation}\label{eq:LHSBLaws-IC}
W(x,0) = W_0(x), \; x \in (0,l),
\end{equation}
for some function $ W_0: (0,l) \rightarrow \mathbb{R}^{k} $ and linear feedback boundary conditions with disturbances defined by 
\begin{equation}\label{eq:LHSBLaws-BCs}
\begin{bmatrix} {W^+}(0,t) \\ {W^-}(l,t) \end{bmatrix} = K\begin{bmatrix} {W^+}(l,t) \\ {W^-}(0,t) \end{bmatrix} + Mb(t), \; t \in (0, +\infty),
\end{equation}
where $  K \in \mathbb{R}^{k\times k}$ is a constant matrix of the form 
$ K = \begin{bmatrix} 0 & K^- \\ K^+ & 0 \end{bmatrix} $, with $ {K^-} \in \mathbb{R}^{m \times (k-m)} $ and $ {K^+} \in \mathbb{R}^{(k-m) \times m} $, $ M \in \mathbb{R}^{k\times k} $ is a non-zero constant diagonal matrix,  and  $ b \in \mathbb{R}^k $ is a vector of disturbance functions. Further more, initial-boundary compatibility conditions are described by
\begin{equation}\label{eq:LHSBLaws-IBCCs}
\begin{bmatrix} {W^+}(0,0) \\ {W^-}(l,0) \end{bmatrix} = K\begin{bmatrix} {W^+}(l,0) \\ {W^-}(0,0) \end{bmatrix}.  
\end{equation}
Note that in the initial-boundary compatibility conditions \eqref{eq:LHSBLaws-IBCCs}, there is no boundary disturbance. That means at the initial time ($ t = 0 $), we assumed there will be no disturbance. It is for such a system that the Input-to-State Stability (ISS) will be discussed in this paper.

In science and engineering, many important physical phenomena, in particular flow of fluids such as flow of shallow water, gas, traffic and electricity, have mathematical models that describe the dynamic behaviour of the flow in terms of mathematical equations. These mathematical models are mainly represented by hyperbolic systems of balance laws, e.g. Saint-Venant equations, isentropic Euler equations, or Telegrapher's equations. The solution of linear hyperbolic systems of balance laws under an initial condition, boundary conditions and initial-boundary compatibility conditions exist and are unique (see \cite{bastin2016stability, prieur2018boundary}). Stabilisation problems with boundary controls (also called boundary feedbacks or boundary damping) of such systems have been an active research field as demonstrated by these papers, \cite{bastin2011boundary, krstic2008backstepping, coron2007strict, de2003boundary, dos2008boundary, diagne2012lyapunov, christofides1996feedback, bastin2008using, coron2015dissipative, gugat2014boundary, gugat2011existence}. These studies mainly focused on linear and non-linear systems in $ L^2- $norm and $ H^2- $norm, respectively,  in the sense of exponential stability. For the most part, a strict Lyapunov function has played a central role in the investigation of conditions  for stability.

Recently, the stabilisation of linear hyperbolic systems of balance laws with boundary disturbance created another dimension in the field. In \cite{lamare2018robust, tanwani2018stabilization}, an input-to-state stability (ISS) which is an exponential stability in the presence of disturbances was introduced for hyperbolic system of conservation laws and balance laws.  

Our aim is to analyse a numerical feedback boundary stabilisation of such systems with boundary disturbance. This method has been presented in a few papers, for instance, \cite{banda2013numerical, dick2014stabilization, gottlich2017numerical, Banda2018, herty2016boundary, gersterdiscretized, gottlich2016electric}. In these studies, a discrete $ L^2- $Lyapunov function is constructed and used to investigate conditions for exponential stability of discretised hyperbolic systems. Furthermore, the decay of the discrete $ L^2- $Lyapunov function has been shown and numerical computations have been done to compare with analytical stability results.

In this paper, we extend our result \cite{Banda2018} in the presence of boundary disturbances. For this reason, we discretise the ISS-Lyapunov function to investigate conditions for ISS in the sense of discrete ISS. Furthermore, the decay of ISS-Lyapunov functions is  explicitly defined. 

This paper is organised as follows: In Section \ref{sec:sec001}, the problem is described. Basic definitions and theoretical results are stated and presented in Section \ref{sec:sec001}. In Section \ref{sec:sec02}, the numerical methods and discretisation are discussed and presented. Also the numerical results are discussed and presented in Section \ref{sec:sec02}. The discussion in Section \ref{sec:sec02} is applied to computational examples in Section \ref{sec:sec03}. Finally, conclusion and references are given at the end.


\section{Preliminaries and analytical results}\label{sec:sec001}

In this section, necessary definitions and theoretical results for the continuous problem will be presented. Firstly, reference will be made to the existence of solutions. This will be followed by a definition of a Lyapunov function and a stability proof in Theorem \ref{thm:Ch01Part2-01}.

In this paper, the sets $ \mathbb{R}^{k} $, $ \mathbb{R}^{k \times k} $ and $ \mathbb{R}_{+}^{k \times k} $ are the set of $ k- $order real vectors, $ k- $order real matrices and $ k- $order positive real matrices, respectively. In addition, the sets $C^0$ and $C^1$ are the set of continuous and once continuously differentiable functions in $ \mathbb{R}^{k} $, respectively. For a given function $ f: [0,l] \rightarrow \mathbb{R}^{k} $, $ L^2- $norm is defined as 
$\displaystyle	\| f \|_{L^2} = \sqrt{\int_{0}^{l} |f(x)|^2dx}$, where $ |\cdot| $ is the Euclidean norm in $ \mathbb{R}^{k} $. Furthermore, $ L^2(0,l) $ is called the space of all measurable functions $ f $ for which $ \| f \|_{L^2} < \infty $.

In order to discuss ISS  of steady-state, $ W \equiv 0 $, of the system \eqref{eq:LHSBLaws} with initial condition \eqref{eq:LHSBLaws-IC}, boundary conditions \eqref{eq:LHSBLaws-BCs} and compatibility conditions \eqref{eq:LHSBLaws-IBCCs}, we make the following assumptions: For all $ x \in [0,l]$, and $ t \in [0,+\infty)$, we assume that 
\begin{enumerate}  
	\item[\textbf{A1.}] The real diagonal matrix $ \Lambda $ is of class $ C^1([0,l]) $.
	\item[\textbf{A2.}] The real matrix $ \Pi $ is of class $ C^0([0,l]) $.
	\item[\textbf{A3.}] The vector of boundary disturbances, $ b $, is a class of $ C^0([0,+\infty)) $.
	\item[\textbf{A4.}] The $\displaystyle \sup_{s \in [0,t]}\left(|b(s)|^2\right) $ is sufficiently small. 
\end{enumerate}

Consider the assumptions \textbf{A1-A4}, existence and uniqueness of a solution to the system \eqref{eq:LHSBLaws} with initial condition \eqref{eq:LHSBLaws-IC}, boundary conditions \eqref{eq:LHSBLaws-BCs} and compatibility conditions \eqref{eq:LHSBLaws-IBCCs} were discussed in detail in \cite{kmit2008classical}. This was accompanied by the proof of existence and uniqueness. For brevity, such details will not be presented in the current paper.

Below, we provide a definition of ISS stability:

\begin{Def}[ISS]\label{def:ISS} 
	The steady-state $ W \equiv 0 $ of the system \eqref{eq:LHSBLaws} with the boundary conditions \eqref{eq:LHSBLaws-BCs} is ISS in  $L^2-$norm with respect to disturbance function $ b $ if there exist positive real constants $\eta > 0$, $ \xi > 0 $, $ C_1 > 0 $ and $ C_2 > 0$ such that, for every initial condition $W_0(x) \in L^2((0,l);\mathbb{R}^k)$  satisfying the compatibility condition \eqref{eq:LHSBLaws-IBCCs}, the  $L^2-$solution to the system \eqref{eq:LHSBLaws} with initial condition \eqref{eq:LHSBLaws-IC}, boundary conditions \eqref{eq:LHSBLaws-BCs} satisfies
	\begin{equation}\label{eq:ExponStabCondChap1-Part2}
	{\|W(\cdot, t) \|}_{L^2((0,l);\mathbb{R}^k)}^2 \leq C_1 {e}^{-\eta t}{\|W_0 \|}_{L^2((0,l);\mathbb{R}^k)}^2 + \frac{C_2}{\eta}\left(1+ \frac{1}{\xi}\right) \sup_{s \in [0,t]}\left(|b(s)|^2\right),\;t \in [0, +\infty).
	\end{equation}	
\end{Def}

\begin{Rk}
\begin{enumerate}
\item The second term on the right hand side (RHS) of the inequality \eqref{eq:ExponStabCondChap1-Part2} estimates the influence of the disturbance function $ b(t) $ on the solution of the system \eqref{eq:LHSBLaws} with the boundary conditions \eqref{eq:LHSBLaws-BCs}.
\item A similar problem was considered in \cite{tanwani2018stabilization} for the case in which $\Lambda$ and $\Pi$ in Equation \eqref{eq:LHSBLaws} are constants.
\item In \cite{tanwani2018stabilization} it was pointed out that stabilisation in the $L^2$-norm does not necessarily guarantee convergence of the maximum norm of $W(\cdot, t)$ over the domain $[0,l]$ in space. To guarantee such convergence, stability is considered in the $H^1$-norm.
\item In this paper, analysis will be made in the $L^2$-norm. 
\end{enumerate} 
\end{Rk}

Similar to Definition \ref{def:ISS}, we define an ISS-Lyapunov function as follows: 
\begin{Def}[$ L^2- $ISS-Lyapunov function]
	For any continuously differentiable weight function defined by $P(x)= \text{diag}\{ {P^+}(x) , {P^-}(x) \} $, where $ P^+(x) \in \mathbb{R}_{+}^{m \times m}$ and $P^-(x) \in \mathbb{R}_{+}^{(k-m) \times (k-m)}$, an $ L^2- $function defined by 
	\begin{equation}\label{eq:LyapunovfunCh01}
	\mathcal{L}(W(\cdot,t)) = \int_{0}^{l} W^{\top} P(x)W dx, \; t \in [0,+\infty),
	\end{equation}  
	is said to be an ISS-Lyapunov function for the system \eqref{eq:LHSBLaws} with the boundary conditions \eqref{eq:LHSBLaws-BCs} if there exist  positive real constants $ \eta > 0 $, $ \xi > 0 $ and $ \nu > 0 $ such that, for all functions $ b(t) \in C^0([0,+\infty))  $, for all solutions of the system \eqref{eq:LHSBLaws} satisfying the boundary conditions \eqref{eq:LHSBLaws-BCs}, and for all $ t \in [0, +\infty) $, 
	\begin{equation}\label{AppxDtLyapunovfun}
	\frac{d\mathcal{L}(W(\cdot,t))}{dt} \leq - \eta \mathcal{L}(W(\cdot,t)) + \nu \left(1 + \frac{1}{\xi}\right) \sup_{s \in [0,t]}\left(|b(s)|^2\right).
	\end{equation}
\end{Def}

The following proposition presents preliminary results which will be used in the proof of the main result of this section in Theorem \ref{thm:Ch01Part2-01}:
 \begin{Prop}\label{prop:Prop1}
	Let $ y$ and $z$ be vectors in $\mathbb{R}^k$. For any matrix $ A $ and positive semi-definite matrix $ B $ in $ \mathbb{R}^{k\times k} $, the following holds:
	\begin{itemize}
		\item[a)] 
		\begin{equation}\label{eq:QuadraticRelation1}
		- 2y^\top A(y - z) = - y^\top A y + z^\top A z - (y -z )^\top A (y -z ).
		\end{equation}
		\item[b)] there exists $ \xi > 0 $ such that 
		\begin{equation}\label{eq:QuadraticRelation2}
		\pm 2y^\top B z \leq \xi y^\top B y + \frac{1}{\xi}z^\top B z.
		\end{equation} 
	\end{itemize}
\end{Prop}
\begin{Pf} The proof of the above statements is straightforward:
	a) Consider a quadratic form to obtain the equation \eqref{eq:QuadraticRelation1} as follows:
	\begin{align*}
	(y -z )^\top A (y -z ) =&\; y^\top A y + z^\top A z - 2y^\top A z,\\
	=&\; - y^\top A y + z^\top A z - 2y^\top A z + 2y^\top A y,\\
	=&\; - y^\top A y + z^\top A z  + 2y^\top A(y - z).
	\end{align*}
	
	b) The following inequality implies the inequality \eqref{eq:QuadraticRelation2}:
	\begin{align*}
	0 \leq&\; \left(\sqrt{\xi}y \mp \frac{1}{\sqrt{\xi}}z\right)^\top B \left(\sqrt{\xi}y \mp \frac{1}{\sqrt{\xi}}z\right),\\
	=&\; {\xi} y^\top B y + \frac{1}{\xi}z^\top B z \mp 2y^\top B z.
	\end{align*}
	\qed	
\end{Pf}

In Lemma \ref{lem:Lemma01-Part2} below, the boundedness of the Lyapunov function is established:

\begin{Lem}\label{lem:Lemma01-Part2}
	Denote the smallest and largest eigenvalues of the diagonal matrix $ P(x) $ by $ \zeta $ and $\beta $, respectively. 
	Then, there exists a positive real constant $ \eta > 0 $, and for every $ W$, the inequalities  
	\begin{align}
	\zeta \int_{0}^{l} |W|^2 dx &\leq \mathcal{L}(W(\cdot,t)) \leq \beta \int_{0}^{l} |W|^2 dx.\label{eq:Lem01-Part2-cond02}
		\end{align}
	hold.
\end{Lem}

\begin{Pf}
Since the diagonal matrix $ P(x)$ is positive definite for all $x \in [0,l] $, for every $ W, $ the following holds: 
\begin{equation}\label{eq:lllllll}
\zeta |W|^2 \leq W^{\top}P(x) W \leq \beta |W|^2,\;  \forall \; W \in \mathbb{R}^k,\; x \in [0,l]. 
\end{equation}
Thus, Inequality \eqref{eq:Lem01-Part2-cond02} is obtained. \qed
\end{Pf}

Further, a version of the well known Gronwall's Lemma is stated as follows: 
\begin{Lem}[Gronwall's Lemma]\label{lem:Lemma02-Part2}
	Let $ y \in C^1([0, +\infty)) $, $ z \in \mathbb{R} $, $ a \in \mathbb{R}^+ $, and 
	\begin{equation*}
	y'(t) \leq - a y(t) + z,\quad y(0) = c \geq 0, \quad t \geq 0.
	\end{equation*}
	Then 
	\begin{equation*}
	y(t) \leq \left(c - \frac{z}{a}\right) e^{-at} + \frac{z}{a}, \quad t \geq 0.
	\end{equation*}
\end{Lem}
\begin{Pf}
The proof of a general case of Gronwall's Lemma is given in  Lemma 1.1.1 in \cite{lakshmikantham1989stability}. Therein the coefficients $ a $ and $ z $ are functions of $ t $. We adopt the proof by considering constants $ a $ and $ z $. \qed
\end{Pf}

We now state the stability result as follows 
\begin{Thm}[Stability]\label{thm:Ch01Part2-01}
	Assume the system \eqref{eq:LHSBLaws} with the boundary conditions \eqref{eq:LHSBLaws-BCs} satisfies assumptions \textbf{A1-A4}. Let $ \xi $ be any positive real number. Assume that the matrix 
	\begin{equation}\label{eq:Theorem01Part2-cond01}
	 -\Lambda(x) P'(x) - \Lambda'(x) P(x) + \Pi^{\top}(x)P(x) + P(x) \Pi(x),
	 \end{equation}
	is positive definite for all $ x \in [0,l]$ where $\Lambda(x)$ are also continuously differentiable and the matrix 
	\begin{equation}\label{eq:Theorem01Part2-cond02}
	\begin{bmatrix} {\Lambda^+}(l){P^+}(l) & 0 \\ 0 &  {\Lambda^-}(0){P^-}(0) \end{bmatrix} - \left(1 + \xi \right)K^{\top}\begin{bmatrix} {\Lambda^+}(0){P^+}(0) & 0 \\ 0 &  {\Lambda^-}(l){P^-}(l) \end{bmatrix}K,
	\end{equation}
	is positive semi-definite. Moreover, let $ \nu $ be the largest eigenvalue of the matrix 
	\begin{equation*}
	M^{\top}\begin{bmatrix} {\Lambda^+}(0){P^+}(0) & 0 \\ 0 &  {\Lambda^-}(l){P^-}(l) \end{bmatrix}M.
	\end{equation*} 
	Then the  $ L^2- $function defined by \eqref{eq:LyapunovfunCh01}
	 is an ISS-Lyapunov function for the system \eqref{eq:LHSBLaws} with boundary conditions \eqref{eq:LHSBLaws-BCs}. Moreover, the steady-state $ W(x,t) \equiv 0 $ of the system \eqref{eq:LHSBLaws} with boundary conditions \eqref{eq:LHSBLaws-BCs} is ISS in  $L^2-$norm with respect to the disturbance function $ b $.
\end{Thm}
\begin{Rk}
There are two commonly used forms of weight functions in the $ L^2- $function \eqref{eq:LyapunovfunCh01}. These are explicit weight functions (see \cite{bastin2017quadratic}) and implicit weight functions (see \cite{bastin2011boundary, coron2007strict, dos2008boundary, diagne2012lyapunov, bastin2008using}). The implicit weight function is defined by 
$ P(x)= \text{diag}\{ {P^+}\exp(-\mu x) , {P^-}\exp(\mu x) \}, \quad \mu > 0$, 
where $ P^+ \in \mathbb{R}_{+}^{m \times m} $ and $ P^- \in \mathbb{R}_{+}^{(k-m) \times (k-m)}$ are constant diagonal matrices.    
\end{Rk}
At this point, we proceed with the proof of Theorem \ref{thm:Ch01Part2-01}
\begin{Pf}
We consider the $ L^2- $function \eqref{eq:LyapunovfunCh01} as a candidate ISS-Lyapunov function. By computing a time derivative of the candidate ISS-Lyapunov function as in \cite{bastin2016stability}(see Section 5.1) and \cite{bastin2008using}, we obtain 
\begin{align}
\frac{d\mathcal{L}(W(\cdot,t))}{dt} =&\;-\left[ W^{\top} \Lambda(x)P(x)W \right]_{0}^{l}\nonumber\\
& \;  -\int_{0}^{l} W^{\top}\left(- \Lambda(x)P'(x) - \Lambda'(x)P(x) + \Pi(x)^{\top}P(x) + P(x)\Pi(x)\right)W dx. \label{eq:Thm1Prof-02}
\end{align}

At this stage the boundary conditions \eqref{eq:LHSBLaws-BCs} and the compatibility conditions \eqref{eq:LHSBLaws-IBCCs} are inserted to obtain:
 {\footnotesize
\begin{align}
-\left[ W^{\top} \Lambda(x)P(x)W \right]_{0}^{l} =&\; -\begin{bmatrix} {W^+}(l,t) \\ {W^-}(0,t) \end{bmatrix}^{\top}\begin{bmatrix} {\Lambda^+}(l){P^+}(l) & 0 \\ 0 &  {\Lambda^-}(0){P^-}(0) \end{bmatrix}\begin{bmatrix} {W^+}(l,t) \\ {W^-}(0,t) \end{bmatrix}\nonumber\\
&\; +\left(K\begin{bmatrix} {W^+}(l,t) \\ {W^-}(0,t) \end{bmatrix} + Mb(t)\right)^{\top}\begin{bmatrix} {\Lambda^+}(0){P^+}(0) & 0 \\ 0 &  {\Lambda^-}(l){P^-}(l) \end{bmatrix}\left(K\begin{bmatrix} {W^+}(l,t) \\ {W^-}(0,t) \end{bmatrix} + Mb(t)\right),\nonumber\\
=&\; -\begin{bmatrix} {W^+}(l,t) \\ {W^-}(0,t) \end{bmatrix}^{\top}\begin{bmatrix} {\Lambda^+}(l){P^+}(l) & 0 \\ 0 &  {\Lambda^-}(0){P^-}(0) \end{bmatrix}\begin{bmatrix} {W^+}(l,t) \\ {W^-}(0,t) \end{bmatrix}\nonumber\\
&\; + \begin{bmatrix} {W^+}(l,t) \\ {W^-}(0,t) \end{bmatrix}^{\top}K^{\top}\begin{bmatrix} {\Lambda^+}(0){P^+}(0) & 0 \\ 0 &  {\Lambda^-}(l){P^-}(l) \end{bmatrix}K\begin{bmatrix} {W^+}(l,t) \\ {W^-}(0,t) \end{bmatrix}\nonumber\\
&\; + 2\left(K\begin{bmatrix} {W^+}(l,t) \\ {W^-}(0,t) \end{bmatrix} \right)^{\top}\begin{bmatrix} {\Lambda^+}(0){P^+}(0) & 0 \\ 0 &  {\Lambda^-}(l){P^-}(l) \end{bmatrix}Mb(t),\nonumber\\
&\; + b(t)^{\top}M^{\top}\begin{bmatrix} {\Lambda^+}(0){P^+}(0) & 0 \\ 0 &  {\Lambda^-}(l){P^-}(l) \end{bmatrix}Mb(t).\label{eq:Thm1Prof-03}
\end{align}
 }

Inequality \eqref{eq:QuadraticRelation2} in Proposition \ref{prop:Prop1} is used for positive semi-definite quadratic form on the RHS of the equation \eqref{eq:Thm1Prof-03} to obtain: 
\begin{align}
-\left[ W^{\top} \Lambda(x)P(x)W \right]_{0}^{l} \leq &\; -\begin{bmatrix} {W^+}(l,t) \\ {W^-}(0,t) \end{bmatrix}^{\top}\begin{bmatrix} {\Lambda^+}(l){P^+}(l) & 0 \\ 0 &  {\Lambda^-}(0){P^-}(0) \end{bmatrix}\begin{bmatrix} {W^+}(l,t) \\ {W^-}(0,t) \end{bmatrix}\nonumber\\
&\; + \left(1 + \xi \right)\begin{bmatrix} {W^+}(l,t) \\ {W^-}(0,t) \end{bmatrix}^{\top}K^{\top}\begin{bmatrix} {\Lambda^+}(0){P^+}(0) & 0 \\ 0 &  {\Lambda^-}(l){P^-}(l) \end{bmatrix}K\begin{bmatrix} {W^+}(l,t) \\ {W^-}(0,t) \end{bmatrix}\nonumber\\
&\; +\left(1 + \frac{1}{\xi}\right) b(t)^{\top}M^{\top}\begin{bmatrix} {\Lambda^+}(0){P^+}(0) & 0 \\ 0 &  {\Lambda^-}(l){P^-}(l) \end{bmatrix}Mb(t),\nonumber\\
= & -\begin{bmatrix} {W^+}(l,t) \\ {W^-}(0,t) \end{bmatrix}^{\top}\left(\begin{bmatrix} {\Lambda^+}(l){P^+}(l) & 0 \\ 0 &  {\Lambda^-}(0){P^-}(0) \end{bmatrix}\right.\nonumber\\
&\; \qquad \qquad \left.- \left(1 + \xi \right)K^{\top}\begin{bmatrix} {\Lambda^+}(0){P^+}(0) & 0 \\ 0 &  {\Lambda^-}(l){P^-}(l) \end{bmatrix}K \right)\begin{bmatrix} {W^+}(l,t) \\ {W^-}(0,t) \end{bmatrix}\nonumber\\
&\; +\left(1 + \frac{1}{\xi}\right) b(t)^{\top}M^{\top}\begin{bmatrix} {\Lambda^+}(0){P^+}(0) & 0 \\ 0 &  {\Lambda^-}(l){P^-}(l) \end{bmatrix}Mb(t).\label{eq:Thm1Prof-04}
\end{align}
Therefore, inserting Equation \eqref{eq:Thm1Prof-04} into Equation \eqref{eq:Thm1Prof-02} gives:
\begin{align}
\frac{d\mathcal{L}(W(\cdot,t))}{dt} \leq &\;  -\int_{0}^{l} W^{\top}\left(- \Lambda(x)P'(x) - \Lambda'(x)P(x) + \Pi(x)^{\top}P(x) + P(x)\Pi(x)\right)W dx \nonumber\\
& -\begin{bmatrix} {W^+}(l,t) \\ {W^-}(0,t) \end{bmatrix}^{\top}\left(\begin{bmatrix} {\Lambda^+}(l){P^+}(l) & 0 \\ 0 &  {\Lambda^-}(0){P^-}(0) \end{bmatrix}\right.\nonumber\\
&\; \qquad \qquad \left.- \left(1 + \xi \right)K^{\top}\begin{bmatrix} {\Lambda^+}(0){P^+}(0) & 0 \\ 0 &  {\Lambda^-}(l){P^-}(l) \end{bmatrix}K \right)\begin{bmatrix} {W^+}(l,t) \\ {W^-}(0,t) \end{bmatrix}\nonumber\\
&\; +\left(1 + \frac{1}{\xi}\right) b(t)^{\top}M^{\top}\begin{bmatrix} {\Lambda^+}(0){P^+}(0) & 0 \\ 0 &  {\Lambda^-}(l){P^-}(l) \end{bmatrix}Mb(t).
 \label{eq:Thm1Prof-05}
\end{align} 

Applying the assumption that $ \nu $ is the largest eigenvalue of the matrix 
	\begin{equation*}
	M^{\top}\begin{bmatrix} {\Lambda^+}(0){P^+}(0) & 0 \\ 0 &  {\Lambda^-}(l){P^-}(l) \end{bmatrix}M,
	\end{equation*} 
using the assumption in Theorem \ref{thm:Ch01Part2-01} for the matrix \eqref{eq:Theorem01Part2-cond02}, Inequality \eqref{eq:Thm1Prof-05} is reduced to 

\begin{align}
\frac{d\mathcal{L}(W(\cdot,t))}{dt} \leq &\;  -\int_{0}^{l} W^{\top}\left(- \Lambda(x)P'(x) - \Lambda'(x)P(x) + \Pi(x)^{\top}P(x) + P(x)\Pi(x)\right)W dx\nonumber \\
&\; + \nu \left(1 + \frac{1}{\xi}\right) |b(t)|^2,\nonumber\\
\leq &\;  -\int_{0}^{l} W^{\top}Q(x)W dx + \nu \left(1 + \frac{1}{\xi}\right)\sup_{s \in [0,t]}\left(|b(s)|^2\right),\label{eq:Thm1Prof-06}
\end{align}	 
where $ Q(x) = - \Lambda(x)P'(x) - \Lambda'(x)P(x) + \Pi(x)^{\top}P(x) + P(x)\Pi(x)$. Furthermore, by the assumption in Theorem \ref{thm:Ch01Part2-01} for the matrix \eqref{eq:Theorem01Part2-cond01}, i.e. positive definiteness of $ Q(x) $, there exist $ \eta > 0 $ such that $ W^{\top}Q(x)W \geq \eta W^{\top}P(x)W $. Thus, the inequality \eqref{eq:Lem01-Part2-cond01} below is obtained:
\begin{align} 
\frac{d\mathcal{L}(W(\cdot,t))}{dt} &\leq - \eta \mathcal{L}(W(\cdot,t)) + \nu \left(1 + \frac{1}{\xi}\right) \sup_{s \in [0,t]}\left(|b(s)|^2\right).\label{eq:Lem01-Part2-cond01}
\end{align}


For the purpose of completing the proof, the Gronwall's Lemma \ref {lem:Lemma02-Part2} is applied to obtain:

\begin{align}
\mathcal{L}(W(\cdot,t)) \leq&\; e^{-\eta t} \left(\mathcal{L}(W(\cdot,0)) - \frac{\nu}{\eta }\left(1+ \frac{1}{\xi}\right)\sup_{s \in [0,t]} \left(|b(s)|^2\right) \right)\nonumber \\
&\; + \frac{\nu}{\eta}\left(1+ \frac{1}{\xi}\right)\sup_{s \in [0,t]} \left(|b(s)|^2\right),\nonumber\\
&\;\leq e^{-\eta t} \mathcal{L}(W(\cdot,0))  + \frac{\nu}{\eta}\left(1+ \frac{1}{\xi}\right)\sup_{s \in [0,t]} \left(|b(s)|^2\right),\;t \geq 0.\label{eq:llllll}
\end{align}

Now insert the inequality in \eqref{eq:Lem01-Part2-cond02} into Inequality \eqref{eq:llllll}, to obtain
\begin{align}
\zeta{\|W(\cdot, t) \|}_{L^2((0,l);\mathbb{R}^k)}^2 \leq \beta {e}^{-\eta t}{\|W_0 \|}_{L^2((0,l);\mathbb{R}^k)}^2  + \frac{\nu}{\eta}\left(1+ \frac{1}{\xi}\right)\sup_{s \in [0,t]} \left(|b(s)|^2\right),\;t \geq 0.\label{eq:llllllll}
\end{align}

Therefore, from the inequality \eqref{eq:llllllll}, the constant coefficients in the condition for exponential stability \eqref{eq:ExponStabCondChap1-Part2} can be assigned to  $ C_1 = \beta/\zeta $ and $ C_2 = \nu/\zeta $, hence Theorem \ref{thm:Ch01Part2-01} is proved.
\qed 
\end{Pf}

Correspondingly, consider a $ k\times k $ uniform linear hyperbolic system of balance laws which can be written as 
\begin{equation}\label{eq:LHSBLaws2}
\partial_t W + \Lambda \partial_x W + \Pi W  = 0,
\end{equation}
where $ \Lambda, M \in \mathbb{R}^{k \times k} $ are non-zero diagonal matrices, and $ \Pi, K \in \mathbb{R}^{k \times k}$ are non-zero matrices with the boundary conditions \eqref{eq:LHSBLaws-BCs}. 

\begin{Cc}\label{col:Col01-Part2}
	Assume the system \eqref{eq:LHSBLaws2} with the boundary conditions \eqref{eq:LHSBLaws-BCs} satisfies assumptions \textbf{A3-A4}. Let $ \xi $ be any positive real number. Assume that the matrix 
	\begin{equation}\label{eq:Col01-cond01}
	-\Lambda P'(x) + \Pi^{\top}P(x) + P(x) \Pi,
	\end{equation}
	is positive definite for all $ x \in [0,l]$, and the matrix  
	\begin{equation}\label{eq:Col01-cond02}
	\begin{bmatrix} \Lambda^+P^+(l) & 0 \\ 0 &  \Lambda^- P^-(0) \end{bmatrix} - \left(1 + \xi \right)K^{\top} \begin{bmatrix} \Lambda^+P^+(0)  & 0 \\ 0 & - \Lambda^- P^-(l) \end{bmatrix} K,
	\end{equation}
	is positive semi-definite. Then the  $ L^2- $function, $ \mathcal{L}, $ defined by \eqref{eq:LyapunovfunCh01} is an ISS-Lyapunov function for the system \eqref{eq:LHSBLaws2} with boundary conditions \eqref{eq:LHSBLaws-BCs}. Moreover, the steady-state $ W(x,t) \equiv 0 $ of the system \eqref{eq:LHSBLaws2} with boundary conditions \eqref{eq:LHSBLaws-BCs} is ISS in  $L^2-$norm with respect to disturbance function $ b $.
\end{Cc}

Having established the stability of the continuous model, Equation \eqref{eq:LHSBLaws}, we now move on to analyse the stability of the discretised form of the same equation in the next section.


\section{Numerical discretisation and stabilisation for a balance law with boundary disturbance}\label{sec:sec02}

The discretisation of the balance law in Equation \eqref{eq:LHSBLaws} will be discussed first. This will be followed by the discrete presentation of the Lyapunov function and the stability analysis of the discrete system. In order to solve a linear hyperbolic system of balance laws numerically,  a time splitting technique which consists of a linear hyperbolic system of conservation laws and a linear system of ordinary differential equations is applied. Thus, the non-uniform system \eqref{eq:LHSBLaws} can be written as follows:
\begin{subequations}\label{eq:LHSBLaws-Split}
	\begin{align}
	\partial_t W + \Lambda (x) \partial_x W =&\;0,\label{eq:LHSBLs-Split-01}\\
	\partial_t W + \Pi(x) W =&\; 0,\label{eq:LHSBLs-Split-02}
	\end{align}
\end{subequations}
where $(x,t) \in [0,l]\times [0,+\infty) $. A first-order Finite Volume Method (\textbf{FVM}), the upwind scheme, is applied to discretise space together with Euler schemes for temporal discretisation. The details of the use of the approach can be found in \cite{leveque2002finite, martinez2018numerical, weldegiyorgis2016numerical}. Specifically, we fix $ T > 0 $ and discretise the domain with $(x,t) \in [0, l]\times [0, T]$ by taking uniform space and time step sizes as $ \Delta x = l/J $ and $ \Delta t = T/N $, where $ J, N > 0 $, respectively. The values $J$ and $N$ denote the number of cells in space and time, respectively. Denote grid points by \[ x_{j - \frac{1}{2}} = j {\Delta x},\; j = 0, \dots, J,\quad  t^n = n {\Delta t},\; n = 0, \dots, N. \] Further, denote left and right boundary points by $ x_{- \frac{1}{2}} = 0 $ and $ x_{J - \frac{1}{2}} = l $, respectively. In addition, cell centres are denoted by $ x_{j} = \left(j + \frac{1}{2}\right) {\Delta x}, \; j = 0, \dots, J-1 $. 

A first order numerical scheme as described in \cite{leveque2002finite} is considered. The approximate cell average of the state variable, $ W $, over the $ j^{\text{th}} $ cell at time $ t^n\;(n =  0, \dots, N) $ is defined by 
\begin{equation}\label{eq:CellAvg}
W_{j}^{n} = \frac{1}{\Delta x}\int_{x_{j - \frac{1}{2}}}^{x_{j + \frac{1}{2}}} W(x,t^n)\; dx, \; j = 0, \dots, J-1,
\end{equation}      
such that for a smooth solution $ W(x,t) $, the integral approximation is defined as  
\begin{equation}\label{eq:IntegralApprox}
\int_{0}^{l} W(x,t^n) dx \approx {\Delta x}\sum_{j = 0}^{J-1}W_{j}^{n}, \;  n = 0, \dots , N-1.
\end{equation}
Therefore, the solution $ W(x_j,t^n) $ is approximated by $ W_j^n $. Hence, for $  n = 0, \dots , N-1 $, $ j = 0, \dots, J-1 $, the non-uniform split system \eqref{eq:LHSBLaws-Split} is discretised as 
\begin{subequations}\label{eq:DiscLHSBLaws}
	\begin{align}
	&\begin{bmatrix} \widetilde{W^+}_{j}^{n}\\\widetilde{W^-}_{j}^{n} \end{bmatrix} =  \begin{bmatrix} {W^+}_{j}^{n}\\{W^-}_{j}^{n} \end{bmatrix}  - {\frac{\Delta t}{\Delta x}} \begin{bmatrix} {\Lambda_{j-1}^+} & 0 \\ 0 & - {\Lambda_{j+1}^-} \end{bmatrix} \begin{bmatrix}  {W^+}_{j}^n - {W^+}_{j-1}^n \\ {W^-}_{j+1}^n - {W^-}_{j}^n \end{bmatrix}, \label{eq:DiscLHSBLaws01}\\
	&\begin{bmatrix} {W^+}_{j}^{n+1}\\{W^-}_{j}^{n+1} \end{bmatrix} =  
	\begin{bmatrix} \widetilde{W^+}_{j}^{n}\\\widetilde{W^-}_{j}^{n} \end{bmatrix} - {\Delta t}\Pi_j \begin{bmatrix} \widetilde{W^+}_{j}^{n}\\\widetilde{W^-}_{j}^{n}  \end{bmatrix}.\label{eq:DiscLHSBLaws02}
	\end{align}
\end{subequations}
Consequently, the initial conditions \eqref{eq:LHSBLaws-IC}, the boundary conditions \eqref{eq:LHSBLaws-BCs} and the compatibility conditions \eqref{eq:LHSBLaws-IBCCs} are discretised as 
\begin{equation}\label{eq:DiscLHSBLaws-IC}
W_j^0 = W_{0,j}, \quad j = 0, \dots,J-1,
\end{equation}
\begin{equation}\label{eq:DiscLHSBLaws-BCs}
\begin{bmatrix} {W^+}_{-1}^{n+1}\\{W^-}_{J}^{n+1} \end{bmatrix} =  K \begin{bmatrix} {W^+}_{J-1}^{n+1}\\{W^-}_{0}^{n+1} \end{bmatrix} + Mb^{n+1},\quad n = 0, \dots , N-1,
\end{equation}
and 
\begin{equation}\label{eq:DiscLHSBLaws-IBCCs}
\begin{bmatrix} {W^+}_{-1}^{0}\\{W^-}_{J}^{0} \end{bmatrix} =  K \begin{bmatrix} {W^+}_{J-1}^{0}\\{W^-}_{0}^{0} \end{bmatrix}, 
\end{equation}
respectively. 

\begin{enumerate}
\item[\textbf{A5.}] Assume that all the assumptions \textbf{A1-A4} hold for the discretised system \eqref{eq:DiscLHSBLaws}.
\end{enumerate}

The aim of this paper is to investigate conditions for numerical boundary feedback stabilisation in the sense of the following definitions of discrete ISS and discrete ISS-Lyapunov function. 

\begin{Def}[Discrete ISS]
	The steady-state $ W_j^n \equiv 0,\; j = 0, \dots, J-1,\; n = 0, \dots, N-1 $ of the discretised system \eqref{eq:DiscLHSBLaws} with the discretised boundary conditions \eqref{eq:DiscLHSBLaws-BCs} is discrete ISS in $L^2-$norm with respect to discrete disturbance function $ b^n,\;n=0,\dots, N-1 $ if there exist positive real constants $\eta > 0$, $ \xi > 0 $, $ C_1 > 0 $ and $C_2 > 0 $ such that, for every initial condition $W_j^0 \in L^2((x_{j-\frac{1}{2}},x_{j+\frac{1}{2}});\mathbb{R}^k),\;j = 0, \dots , J-1$ satisfying the compatibility condition \eqref{eq:DiscLHSBLaws-IBCCs}, the  $L^2-$solution of the discretised system \eqref{eq:DiscLHSBLaws} with initial condition \eqref{eq:DiscLHSBLaws-IC} and boundary conditions \eqref{eq:DiscLHSBLaws-BCs} satisfies
	\begin{equation}\label{eq:DiscExponStabCondChap1-Part2}
	{\Delta x}{\sum_{j=0}^{J-1}} |W_j^{n+1}|^2 \leq C_1e^{-\eta t^{n+1}}{\Delta x}{\sum_{j=0}^{J-1}} |W_j^0|^2 + \frac{C_2}{\eta}\left(1 + \frac{1}{\xi}\right) \sup_{0 \leq s \leq n} \left(|b^s|^2\right), \; n = 0, \dots , N-1.
	\end{equation}
\end{Def}

\begin{Def}[A discrete $ L^2 -$ISS-Lyapunov function]
For any discrete weight function defined by $P_j= \text{diag}\{ P_j^+ , P_j^-\}$, where $ P_j^+$ and $P_j^-$ denote the first $ m $ and the last $ k-m $  positive diagonal entries, respectively, for $j = 0, \dots, J-1$, a discrete $ L^2- $function defined by 
	\begin{equation}\label{eq:DiscLyapunovfunCh01}
	\mathcal{L}^n = {\Delta x }\sum_{j=0}^{J-1} {W_j^n}^{\top} P_j W_j^n, \quad n = 0, \dots , N-1,
	\end{equation} 
	is said to be a discrete ISS-Lyapunov function for the discretised system \eqref{eq:DiscLHSBLaws} with the discretised boundary conditions \eqref{eq:DiscLHSBLaws-BCs} if there exist positive real constants $ \eta > 0 $, $ \xi > 0 $  and $ \nu > 0 $ such that, for all discrete functions  $ b^n,\;n=0,\dots, N-1 $,  for all solutions of the discretised system \eqref{eq:DiscLHSBLaws} satisfying the discretised boundary conditions \eqref{eq:DiscLHSBLaws-BCs}, and for all $ n = 0, \dots, N-1 $, 
	\begin{equation}\label{AppxDtDiscLyapunovfun-Part2}
	\frac{\mathcal{L}^{n+1} - \mathcal{L}^{n}}{\Delta t} \leq - \eta \mathcal{L}^{n} + \nu \left(1 + \frac{1}{\xi}\right) \sup_{0 \leq s \leq n} \left(|b^s|^2\right).
	\end{equation}
\end{Def}

Before stating the main theorem of this section, we present two preliminary results:

\begin{Lem}\label{lem:DiscLemma1}
Assume that $ P_j,\ j = 0, \dots, J-1$, is a positive definite matrix. Define the smallest and largest eigenvalue of $ P_j$ by $\displaystyle \zeta  = \min_{0 \leq j \leq J-1}P_j $ and $\displaystyle \beta  = \max_{0 \leq j \leq J-1}P_j $, respectively. 
Then, the following inequality holds:  
\begin{align}
\zeta {\Delta x }\sum_{j=0}^{J-1} |W_j^n|^2  &\leq \mathcal{L}^{n} \leq \beta  {\Delta x }\sum_{j=0}^{J-1} |W_j^n|^2,\label{eq:DiscLem01-cond02}
\end{align}
\end{Lem}
\begin{Pf}
Since the diagonal matrix $ P_j,\ j = 0, \dots, J-1$ is positive definite, for all $ W_j^n,\ n = 0, \dots, N-1 $, we have 
\begin{equation}\label{eq:PostiveSD02}
\zeta |W_j^n|^2  \leq {W}_{j}^{n \top} P_j {W}_{j}^{n} \leq \beta |W_j^n|^2,\; j = 0, \dots, J-1.
\end{equation}
Then, the inequality \eqref{eq:PostiveSD02} implies the inequality \eqref{eq:DiscLem01-cond02}. 
\qed
\end{Pf}
Now we present an equivalent Gronwall's Lemma for the discrete case:
\begin{Lem}\label{lem:DiscLemma2}
Let $ a > 0 $ and $ z \in \mathbb{R} $. Suppose for discrete functions $ y^{n}, \; n = 0 , \dots, N-1  $, 
\begin{equation}\label{eq:DiscDInequality}
\frac{y^{n+1} - y^{n}}{\Delta t} \leq - a y^{n} + z,\quad y^{0} = c. 
\end{equation}
Then
\begin{equation}\label{eq:DiscInequality}
y^{n+1} \leq \left(c - \frac{z}{a}\right) \left(1 - a{\Delta t} \right)^{n+1} + \frac{z}{a}, \quad n = 0 , \dots, N-1.
\end{equation}
\end{Lem}

\begin{Pf}
	By recursively applying the inequality \eqref{eq:DiscDInequality}, we obtain
	\begin{equation}\label{eq:DiscDInequality1}
	y^{n+1} \leq c \left(1 - a{\Delta t} \right)^{n+1} + z{\Delta t} \sum_{r = 0}^{n} \left(1 - a{\Delta t} \right)^{r}, \; \; n = 0 , \dots, N-1.
	\end{equation}
	Then, the inequality \eqref{eq:DiscDInequality1} implies the inequality \eqref{eq:DiscInequality} for sufficiently small $ {\Delta t}$, $ 0 < 1 - a{\Delta t} < 1 $.
	\qed
\end{Pf}

In the sense of the definitions of discrete ISS and discrete $ L^2 -$ISS-Lyapunov function, we state the numerical stability result as follows: 

\begin{Thm}[Stability]\label{thm:ThmChap01-02-Part2}
Assume the discretised system \eqref{eq:DiscLHSBLaws} with the discretised boundary conditions \eqref{eq:DiscLHSBLaws-BCs} satisfies assumption \textbf{A5}. Let $ T > 0 $ be fixed and the CFL condition, $ \frac{\Delta t}{\Delta x} \max_{\stackrel{1\leq i \leq k}{0 \leq j \leq J-1}} |\lambda_{i,j}| \leq 1 $ hold. Let $ \xi $ be any positive real number. Assume that the matrix
\begin{equation}\label{eq:DiscTheorem01-cond00}
\begin{bmatrix} -{\Lambda_{j-1}^+}\left(\frac{P_{j+1}^+ - P_j^+}{\Delta x}\right) - \left(\frac{{\Lambda_{j}^+} - {\Lambda_{j-1}^+}}{\Delta x}\right)P_{j+1}^+ & 0 \\ 0 & {\Lambda_{j+1}^-}\left(\frac{P_j^- - P_{j-1}^-}{\Delta x}\right) + \left(\frac{{\Lambda_{j+1}^-} -{\Lambda_{j}^-}}{\Delta x}\right)P_{j-1}^- \end{bmatrix},
\end{equation}
is positive definite for all $ j = 0, \dots, J-1 $, and the matrices 
\begin{equation}\label{eq:DiscTheorem01-cond01}
P_j\Pi_j + {\Pi_j}^{\top} P_j - {\Delta t}{\Pi_j}^{\top}P_j{\Pi_j},
\end{equation}
and 
\begin{equation}\label{eq:DiscTheorem01-cond02}
\begin{bmatrix} {\Lambda_{J-1}^+}P_{J}^+ & 0 \\ 0 &  {\Lambda_{0}^-}P_{-1}^- \end{bmatrix} - \left(1 + {\xi}\right)K^{\top}\begin{bmatrix} {\Lambda_{-1}^+}P_{0}^+ & 0 \\ 0 &  {\Lambda_{J}^-}P_{J-1}^- \end{bmatrix} K,
\end{equation}
are positive semi-definite for all $ j = 0, \dots, J-1 $. Then the discrete $ L^2-$function defined by \eqref{eq:DiscLyapunovfunCh01} 
is a discrete ISS-Lyapunov function for the discretised system \eqref{eq:DiscLHSBLaws} with discretised boundary conditions \eqref{eq:DiscLHSBLaws-BCs}. Moreover, the steady-state $ W_j^n \equiv 0,\; j = 0, \dots, J-1,\; n = 0, \dots, N-1 $ of the discretised system \eqref{eq:DiscLHSBLaws} with discretised boundary conditions \eqref{eq:DiscLHSBLaws-BCs} is discrete ISS in $L^2-$norm with respect to discrete disturbance function $ b^n,\;n=0,\dots, N-1 $.
\end{Thm}
\begin{Pf}
The discrete $ L^2-$function \eqref{eq:DiscLyapunovfunCh01} is used to approximate the time derivative of the candidate ISS-Lyapunov $ L^2-$function \eqref{eq:LyapunovfunCh01}. Corresponding to the discrete split system \eqref{eq:DiscLHSBLaws}, the time derivative is approximated in a split form as computed in \cite{Banda2018}. Thus, 
\begin{align}
\frac{\mathcal{L}^{n+1} - \mathcal{L}^{n}}{\Delta t} \leq& - {\Delta x} \sum_{j=0}^{J-1} \widetilde{W}_j^{n \top} \left(P_j\Pi_j + {\Pi_j}^{\top} P_j - {\Delta t}{\Pi_j}^{\top} P_j {\Pi_j} \right) \widetilde{W}_j^n \nonumber \\
& - \sum_{j=0}^{J-1} \begin{bmatrix} {W^+}_{j}^{n}\\{W^-}_{j}^{n} \end{bmatrix}^{\top}  \begin{bmatrix} {\Lambda_{j-1}^+}P_j^+ & 0 \\ 0 & {\Lambda_{j+1}^-}P_j^- \end{bmatrix} \begin{bmatrix}  {W^+}_{j}^n \\ {W^-}_{j}^n \end{bmatrix}\nonumber\\
&\; + \sum_{j=0}^{J-1} \begin{bmatrix} {W^+}_{j-1}^{n}\\{W^-}_{j+1}^{n} \end{bmatrix}^{\top}  \begin{bmatrix} {\Lambda_{j-1}^+}P_j^+ & 0 \\ 0 &  {\Lambda_{j+1}^-}P_j^- \end{bmatrix} \begin{bmatrix}  {W^+}_{j-1}^n \\ {W^-}_{j+1}^n \end{bmatrix},\; n = 0, \dots , N-1.\label{eq:DtDiscLyapunovfun}
\end{align}

By using $ x_j = x_{j-1} + {\Delta x}, \; j= 0,\dots,J-1 $, we obtain, for all $ n = 0, \dots , N-1$: 
\begin{align}
\sum_{j=0}^{J-1} \begin{bmatrix} {W^+}_{j-1}^{n}\\{W^-}_{j+1}^{n} \end{bmatrix}^{\top}  \begin{bmatrix} {\Lambda_{j-1}^+}P_j^+ & 0 \\ 0 &  {\Lambda_{j+1}^-}P_j^- \end{bmatrix} \begin{bmatrix}  {W^+}_{j-1}^n \\ {W^-}_{j+1}^n \end{bmatrix} =&\; \sum_{j=0}^{J-1} \begin{bmatrix} {W^+}_{j}^{n}\\{W^-}_{j}^{n} \end{bmatrix}^{\top}  \begin{bmatrix} {\Lambda_{j}^+}P_{j+1}^+ & 0 \\ 0 &  {\Lambda_{j}^-}P_{j-1}^- \end{bmatrix} \begin{bmatrix}  {W^+}_{j}^n \\ {W^-}_{j}^n \end{bmatrix}\nonumber\\
&\; - \begin{bmatrix} {W^+}_{J-1}^{n}\\{W^-}_{0}^{n} \end{bmatrix}^{\top}  \begin{bmatrix} {\Lambda_{J-1}^+}P_{J}^+ & 0 \\ 0 &  {\Lambda_{0}^-}P_{-1}^- \end{bmatrix} \begin{bmatrix}  {W^+}_{J-1}^n \\ {W^-}_{0}^n \end{bmatrix}\nonumber\\
&\;+ \begin{bmatrix} {W^+}_{-1}^{n}\\{W^-}_{J}^{n} \end{bmatrix}^{\top}  \begin{bmatrix} {\Lambda_{-1}^+}P_{0}^+ & 0 \\ 0 &  {\Lambda_{J}^-}P_{J-1}^- \end{bmatrix} \begin{bmatrix}  {W^+}_{-1}^n \\ {W^-}_{J}^n \end{bmatrix}.\label{eq:DiscThm1Prof-06}
\end{align}
Then, Equation \eqref{eq:DiscThm1Prof-06} is substituted into the Inequality \eqref{eq:DtDiscLyapunovfun} to obtain: 
\begin{align}
\frac{\mathcal{L}^{n+1} - \mathcal{L}^{n}}{\Delta t} \leq& - {\Delta x} \sum_{j=0}^{J-1} \widetilde{W}_j^{n \top} \left(P_j\Pi_j + {\Pi_j}^{\top} P_j - {\Delta t}{\Pi_j}^{\top} P_j {\Pi_j} \right) \widetilde{W}_j^n \nonumber \\
&\; 
- \sum_{j=0}^{J-1} \begin{bmatrix} {W^+}_{j}^{n}\\{W^-}_{j}^{n} \end{bmatrix}^{\top}  \begin{bmatrix} {\Lambda_{j-1}^+}P_j^+ & 0 \\ 0 & {\Lambda_{j+1}^-}P_j^- \end{bmatrix} \begin{bmatrix}  {W^+}_{j}^n \\ {W^-}_{j}^n \end{bmatrix}\nonumber\\
&\; + \sum_{j=0}^{J-1} \begin{bmatrix} {W^+}_{j}^{n}\\{W^-}_{j}^{n} \end{bmatrix}^{\top}  \begin{bmatrix} {\Lambda_{j}^+}P_{j+1}^+ & 0 \\ 0 &  {\Lambda_{j}^-}P_{j-1}^- \end{bmatrix} \begin{bmatrix}  {W^+}_{j}^n \\ {W^-}_{j}^n \end{bmatrix}\nonumber\\
&\; - \begin{bmatrix} {W^+}_{J-1}^{n}\\{W^-}_{0}^{n} \end{bmatrix}^{\top}  \begin{bmatrix} {\Lambda_{J-1}^+}P_{J}^+ & 0 \\ 0 &  {\Lambda_{0}^-}P_{-1}^- \end{bmatrix} \begin{bmatrix}  {W^+}_{J-1}^n \\ {W^-}_{0}^n \end{bmatrix}\nonumber\\
&\;+ \begin{bmatrix} {W^+}_{-1}^{n}\\{W^-}_{J}^{n} \end{bmatrix}^{\top}  \begin{bmatrix} {\Lambda_{-1}^+}P_{0}^+ & 0 \\ 0 &  {\Lambda_{J}^-}P_{J-1}^- \end{bmatrix} \begin{bmatrix}  {W^+}_{-1}^n \\ {W^-}_{J}^n \end{bmatrix}.\label{eq:DiscThm1Prof-07}
\end{align}
for all $ n = 0, \dots , N-1$, The boundary conditions \eqref{eq:DiscLHSBLaws-BCs}, the compatibility conditions \eqref{eq:DiscLHSBLaws-IBCCs}, the inequality \eqref{eq:QuadraticRelation2} in Proposition \ref{prop:Prop1} and the assumption in Theorem \eqref{thm:ThmChap01-02-Part2} are used to simplify the boundary term in the inequality \eqref{eq:DiscThm1Prof-07} as follows:
\begin{align}
&\; - \begin{bmatrix} {W^+}_{J-1}^{n}\\{W^-}_{0}^{n} \end{bmatrix}^{\top}  \begin{bmatrix} {\Lambda_{J-1}^+}P_{J}^+ & 0 \\ 0 &  {\Lambda_{0}^-}P_{-1}^- \end{bmatrix} \begin{bmatrix}  {W^+}_{J-1}^n \\ {W^-}_{0}^n \end{bmatrix} + \begin{bmatrix} {W^+}_{-1}^{n}\\{W^-}_{J}^{n} \end{bmatrix}^{\top}  \begin{bmatrix} {\Lambda_{-1}^+}P_{0}^+ & 0 \\ 0 &  {\Lambda_{J}^-}P_{J-1}^- \end{bmatrix} \begin{bmatrix}  {W^+}_{-1}^n \\ {W^-}_{J}^n \end{bmatrix}  \nonumber\\
=&\; - \begin{bmatrix} {W^+}_{J-1}^{n}\\{W^-}_{0}^{n} \end{bmatrix}^{\top}  \begin{bmatrix} {\Lambda_{J-1}^+}P_{J}^+ & 0 \\ 0 &  {\Lambda_{0}^-}P_{-1}^- \end{bmatrix} \begin{bmatrix}  {W^+}_{J-1}^n \\ {W^-}_{0}^n \end{bmatrix}\nonumber\\
&\; + \left( K \begin{bmatrix} {W^+}_{J-1}^{n}\\{W^-}_{0}^{n} \end{bmatrix} + Mb^n \right)^{\top}\begin{bmatrix} {\Lambda_{-1}^+}P_{0}^+ & 0 \\ 0 &  {\Lambda_{J}^-}P_{J-1}^- \end{bmatrix}\left( K \begin{bmatrix} {W^+}_{J-1}^{n}\\{W^-}_{0}^{n} \end{bmatrix} + Mb^n \right) \nonumber\\
= &\; - \begin{bmatrix} {W^+}_{J-1}^{n}\\{W^-}_{0}^{n} \end{bmatrix}^{\top}  \begin{bmatrix} {\Lambda_{J-1}^+}P_{J}^+ & 0 \\ 0 &  {\Lambda_{0}^-}P_{-1}^- \end{bmatrix} \begin{bmatrix}  {W^+}_{J-1}^n \\ {W^-}_{0}^n \end{bmatrix} + \begin{bmatrix} {W^+}_{J-1}^{n}\\{W^-}_{0}^{n} \end{bmatrix}^{\top}K^{\top}\begin{bmatrix} {\Lambda_{-1}^+}P_{0}^+ & 0 \\ 0 &  {\Lambda_{J}^-}P_{J-1}^- \end{bmatrix} K \begin{bmatrix} {W^+}_{J-1}^{n}\\{W^-}_{0}^{n} \end{bmatrix}\nonumber\\
&\; + 2\left( K \begin{bmatrix} {W^+}_{J-1}^{n}\\{W^-}_{0}^{n} \end{bmatrix} \right)^{\top}\begin{bmatrix} {\Lambda_{-1}^+}P_{0}^+ & 0 \\ 0 &  {\Lambda_{J}^-}P_{J-1}^- \end{bmatrix}\left(Mb^n \right) + b^{n \top} M^{\top} \begin{bmatrix} {\Lambda_{-1}^+}P_{0}^+ & 0 \\ 0 &  {\Lambda_{J}^-}P_{J-1}^- \end{bmatrix} Mb^n,\nonumber\\
\leq &\; - \begin{bmatrix} {W^+}_{J-1}^{n}\\{W^-}_{0}^{n} \end{bmatrix}^{\top}  \begin{bmatrix} {\Lambda_{J-1}^+}P_{J}^+ & 0 \\ 0 &  {\Lambda_{0}^-}P_{-1}^- \end{bmatrix} \begin{bmatrix}  {W^+}_{J-1}^n \\ {W^-}_{0}^n \end{bmatrix} + \begin{bmatrix} {W^+}_{J-1}^{n}\\{W^-}_{0}^{n} \end{bmatrix}^{\top}K^{\top}\begin{bmatrix} {\Lambda_{-1}^+}P_{0}^+ & 0 \\ 0 &  {\Lambda_{J}^-}P_{J-1}^- \end{bmatrix} K \begin{bmatrix} {W^+}_{J-1}^{n}\\{W^-}_{0}^{n} \end{bmatrix}\nonumber\\
&\; + {\xi} \begin{bmatrix} {W^+}_{J-1}^{n}\\{W^-}_{0}^{n} \end{bmatrix}^{\top}K^{\top}\begin{bmatrix} {\Lambda_{-1}^+}P_{0}^+ & 0 \\ 0 &  {\Lambda_{J}^-}P_{J-1}^- \end{bmatrix} K \begin{bmatrix} {W^+}_{J-1}^{n}\\{W^-}_{0}^{n} \end{bmatrix}\nonumber\\
&\; + \frac{1}{\xi}b^{n \top} M^{\top}\begin{bmatrix} {\Lambda_{-1}^+}P_{0}^+ & 0 \\ 0 &  {\Lambda_{J}^-}P_{J-1}^- \end{bmatrix}Mb^n + b^{n \top} M^{\top} \begin{bmatrix} {\Lambda_{-1}^+}P_{0}^+ & 0 \\ 0 &  {\Lambda_{J}^-}P_{J-1}^- \end{bmatrix} Mb^n,\nonumber\\
=&\; - \begin{bmatrix} {W^+}_{J-1}^{n}\\{W^-}_{0}^{n} \end{bmatrix}^{\top}  \begin{bmatrix} {\Lambda_{J-1}^+}P_{J}^+ & 0 \\ 0 &  {\Lambda_{0}^-}P_{-1}^- \end{bmatrix} \begin{bmatrix}  {W^+}_{J-1}^n \\ {W^-}_{0}^n \end{bmatrix}\nonumber\\
&\; + \left(1 + {\xi}\right) \begin{bmatrix} {W^+}_{J-1}^{n}\\{W^-}_{0}^{n} \end{bmatrix}^{\top}K^{\top}\begin{bmatrix} {\Lambda_{-1}^+}P_{0}^+ & 0 \\ 0 &  {\Lambda_{J}^-}P_{J-1}^- \end{bmatrix} K \begin{bmatrix} {W^+}_{J-1}^{n}\\{W^-}_{0}^{n} \end{bmatrix}\nonumber\\
&\; + \left(1 + \frac{1}{\xi}\right)b^{n \top} M^{\top}\begin{bmatrix} {\Lambda_{-1}^+}P_{0}^+ & 0 \\ 0 &  {\Lambda_{J}^-}P_{J-1}^- \end{bmatrix}Mb^n,\nonumber\\
\leq &\; \left(1 + \frac{1}{\xi}\right)b^{n \top} M^{\top}\begin{bmatrix} {\Lambda_{-1}^+}P_{0}^+ & 0 \\ 0 &  {\Lambda_{J}^-}P_{J-1}^- \end{bmatrix}Mb^n,\;n = 0, \dots , N-1. \label{eq:DiscThm1Prof-08}
\end{align}
Thus, applying inequality \eqref{eq:DiscThm1Prof-08}, for all $ n = 0, \dots , N-1 $, inequality \eqref{eq:DiscThm1Prof-07} is simplified as:
\begin{align}
\frac{\mathcal{L}^{n+1} - \mathcal{L}^{n}}{\Delta t} \leq&\; - {\Delta x} \sum_{j=0}^{J-1} \widetilde{W}_j^{n \top} \left(P_j\Pi_j + {\Pi_j}^{\top} P_j - {\Delta t}{\Pi_j}^{\top} P_j {\Pi_j} \right) \widetilde{W}_j^n \nonumber \\
&\; 
- \sum_{j=0}^{J-1} \begin{bmatrix} {W^+}_{j}^{n}\\{W^-}_{j}^{n} \end{bmatrix}^{\top}  \begin{bmatrix} {\Lambda_{j-1}^+}P_j^+ - {\Lambda_{j}^+}P_{j+1}^+ & 0 \\ 0 & {\Lambda_{j+1}^-}P_j^- -{\Lambda_{j}^-}P_{j-1}^- \end{bmatrix} \begin{bmatrix}  {W^+}_{j}^n \\ {W^-}_{j}^n \end{bmatrix}\nonumber\\
&\; + \left(1 + \frac{1}{\xi}\right)b^{n \top} M^{\top}\begin{bmatrix} {\Lambda_{-1}^+}P_{0}^+ & 0 \\ 0 &  {\Lambda_{J}^-}P_{J-1}^- \end{bmatrix}Mb^n,\nonumber\\
= &\; - {\Delta x} \sum_{j=0}^{J-1} \widetilde{W}_j^{n \top} \left(P_j\Pi_j + {\Pi_j}^{\top} P_j - {\Delta t}{\Pi_j}^{\top} P_j {\Pi_j} \right) \widetilde{W}_j^n \nonumber \\
&\; - {\Delta x}\sum_{j=0}^{J-1} W_j^{n \top} {\Theta}_j W_j^n + \left(1 + \frac{1}{\xi}\right)b^{n \top} M^{\top}\begin{bmatrix} {\Lambda_{-1}^+}P_{0}^+ & 0 \\ 0 &  {\Lambda_{J}^-}P_{J-1}^- \end{bmatrix}Mb^n, \label{eq:DiscThm1Prof-09}
\end{align}
where for $ j = 0, \dots, J-1 $,
\[ {\Theta}_j = \begin{bmatrix} -{\Lambda_{j-1}^+}\left(\frac{P_{j+1}^+ - P_j^+}{\Delta x}\right) - \left(\frac{{\Lambda_{j}^+} - {\Lambda_{j-1}^+}}{\Delta x}\right)P_{j+1}^+ & 0 \\ 0 & {\Lambda_{j+1}^-}\left(\frac{P_j^- - P_{j-1}^-}{\Delta x}\right) + \left(\frac{{\Lambda_{j+1}^-} -{\Lambda_{j}^-}}{\Delta x}\right)P_{j-1}^- \end{bmatrix}.\]

Using the result in Lemma \ref{lem:DiscLemma1}, let $ \nu $ be the largest eigenvalue of the matrix  
\begin{equation}\label{eq:boundarymatrix}
M^{\top}\begin{bmatrix} {\Lambda_{-1}^+}P_{0}^+ & 0 \\ 0 &  {\Lambda_{J}^-}P_{J-1}^- \end{bmatrix}M.
 \end{equation}
Furthermore, there exist a positive real number $ \eta > 0 $ (it is explicitly defined in Section \ref{sec:sec03}), by assumption above, such that for every $ W_j^n,\ n = 0, \dots, N-1 $, we have  $ W_{j}^{n \top}{\Theta}_j  W_j^n \geq {\eta} W_{j}^{n \top}P_jW_j^n $. In addition, for $ n = 0, \dots, N-1 $, we have
\begin{equation*}
b^{n \top} M^{\top}\begin{bmatrix} {\Lambda_{-1}^+}P_{0}^+ & 0 \\ 0 &  {\Lambda_{J}^-}P_{J-1}^- \end{bmatrix}Mb^n \leq {\nu} |b^n|^2 \leq {\nu}\sup_{0 \leq s \leq n}\left( |b^s|^2\right). 
\end{equation*} 
Hence, the inequality \eqref{eq:DiscThm1Prof-09} is approximated as 
\begin{equation}\label{eq:DtDiscLyapunovfun00}
\frac{\mathcal{L}^{n+1} - \mathcal{L}^{n}}{\Delta t} \leq - {\eta}\mathcal{L}^{n} + {\nu}\left(1 + \frac{1}{\xi}\right) \sup_{0 \leq s \leq n}\left( |b^s|^2\right),\; n = 0, \dots, N-1.
\end{equation}

From Lemma \ref{lem:DiscLemma2} and by using $\left(1 - \eta {\Delta t} \right)^{n+1} \leq e^{-{\eta}{t}^{n+1}} $, we have 
\begin{align}
\mathcal{L}^{n+1} \leq &\; \left(\mathcal{L}^{0} - \frac{\nu}{\eta}\left(1 + \frac{1}{\xi}\right) \sup_{0 \leq s \leq n}\left(|b^s|^2\right)\right)\left(1 - \eta {\Delta t} \right)^{n+1}\nonumber \\
&\;+ \frac{\nu}{\eta}\left(1 + \frac{1}{\xi}\right) \sup_{0 \leq s \leq n}\left( |b^s|^2\right),\nonumber\\
 \leq &\; e^{-{\eta}{t}^{n+1}}\mathcal{L}^{0} + \frac{\nu}{\eta}\left(1 + \frac{1}{\xi}\right) \sup_{0 \leq s \leq n}\left(|b^s|^2\right), \; n = 0, \dots, N-1.\label{eq:DtDiscLyapunovfun-04}
\end{align}

Thus, for all $ j = 0, \dots, J-1 $, and $ n = 0, \dots, N-1 $, by using the inequalities  \eqref{eq:DiscLem01-cond02} and \eqref{eq:DtDiscLyapunovfun-04}, we have
\begin{equation}\label{eq:DtDiscLyapunovfun-05}
\zeta {\Delta x }\sum_{j=0}^{J-1} |W_j^{n+1}|^2 \leq \beta e^{-{\eta}{t}^{n+1}} {\Delta x }\sum_{j=0}^{J-1} |W_j^0|^2 + \frac{\nu}{\eta}\left(1 + \frac{1}{\xi}\right) \sup_{0 \leq s \leq n}\left(|b^s|^2\right).  
\end{equation}

Therefore, to show that the inequality \eqref{eq:DtDiscLyapunovfun-05} implies the condition for the discrete ISS \eqref{eq:DiscExponStabCondChap1-Part2}, we let $ C_1 = \beta/\zeta $ and $ C_2 = \nu/\zeta $. Hence, the proof of Theorem \ref{thm:ThmChap01-02-Part2} is completed. \qed 
\end{Pf}

Similar to the split system in Equation \eqref{eq:LHSBLaws-Split}, the uniform system \eqref{eq:LHSBLaws2} can also be split as
\begin{subequations}\label{eq:LHSBLaws2-Split}
	\begin{align}
	\partial_t W + \Lambda \partial_x W =&\;0,\label{eq:LHSBLs2-Split-01}\\
	\partial_t W + \Pi W =&\; 0,\label{eq:LHSBLs-Split2-02}
	\end{align}
\end{subequations}
where $(x,t) \in [0,l]\times [0,+\infty) $. Then, the split system \eqref{eq:LHSBLaws2-Split} is discretised as follows 
\begin{subequations}\label{eq:DiscLHSBLaws2}
	\begin{align}
	&\begin{bmatrix} \widetilde{W^+}_{j}^{n}\\\widetilde{W^-}_{j}^{n} \end{bmatrix} =  \begin{bmatrix} {W^+}_{j}^{n}\\{W^-}_{j}^{n} \end{bmatrix}  - {\frac{\Delta t}{\Delta x}} \begin{bmatrix} {\Lambda^+} & 0 \\ 0 & - {\Lambda^-} \end{bmatrix} \begin{bmatrix}  {W^+}_{j}^n - {W^+}_{j-1}^n \\ {W^-}_{j+1}^n - {W^-}_{j}^n \end{bmatrix}, \label{eq:DiscLHSBLaws01cc}\\
	&\begin{bmatrix} {W^+}_{j}^{n+1}\\{W^-}_{j}^{n+1} \end{bmatrix} =  
	\begin{bmatrix} \widetilde{W^+}_{j}^{n}\\\widetilde{W^-}_{j}^{n} \end{bmatrix} - {\Delta t}\Pi \begin{bmatrix} \widetilde{W^+}_{j}^{n}\\\widetilde{W^-}_{j}^{n}  \end{bmatrix}.\label{eq:DiscLHSBLaws02cc}
	\end{align}
\end{subequations}

\begin{Cc}\label{col:Col02-Part2}
	Assume the system \eqref{eq:DiscLHSBLaws2} with boundary conditions \eqref{eq:DiscLHSBLaws-BCs} satisfies {{assumption \textbf{A5}.}} Let $ T > 0 $ be fixed and the CFL condition, $ \frac{\Delta t}{\Delta x} \max_{1\leq i \leq k} |\lambda_{i}| \leq 1 $ hold. Further, let $ \xi $ be any positive real number. Assume that the matrix 
	\begin{equation}\label{eq:DiscCol01-cond00}
	\frac{1}{\Delta x}\begin{bmatrix} -{\Lambda^+}\left(P_{j+1}^+ - P_j^+\right) & 0 \\ 0 & {\Lambda^-}\left(P_j^- - P_{j-1}^-\right) \end{bmatrix},
	\end{equation}
	is positive definite and the matrices
	\begin{equation}\label{eq:DiscCol01-cond01}
	P_j\Pi  + {\Pi}^{\top} P_j - {\Delta t}{\Pi_j}^{\top}P_j{\Pi_j},
	\end{equation}
	and   
	\begin{equation}\label{eq:DiscCol01-cond02}
	\begin{bmatrix} {\Lambda^+}P_{J}^+ & 0 \\ 0 &  {\Lambda^-}P_{-1}^- \end{bmatrix} - \left(1 + {\xi}\right)K^{\top}\begin{bmatrix} {\Lambda^+}P_{0}^+ & 0 \\ 0 &  {\Lambda^-}P_{J-1}^- \end{bmatrix} K,
	\end{equation}
	are positive semi-definite for all $ j = 0, \dots, J-1$. Then the discrete $ L^2-$function defined by \eqref{eq:DiscLyapunovfunCh01} is a discrete ISS-Lyapunov function for the system \eqref{eq:DiscLHSBLaws2} with boundary conditions \eqref{eq:DiscLHSBLaws-BCs}. Moreover, system \eqref{eq:DiscLHSBLaws2} with boundary conditions \eqref{eq:DiscLHSBLaws-BCs} is discrete ISS in $L^2-$norm with respect to discrete disturbance function $ b^n,\;n=0,\dots, N-1 $.
\end{Cc}
{{The proof of Corollary \ref{col:Col02-Part2} is a special case of the proof of Theorem \ref{thm:ThmChap01-02-Part2} for the discretised system \eqref{eq:DiscLHSBLaws2}.}}  The case in Equation \eqref{eq:LHSBLaws2} above was analysed in \cite{tanwani2018stabilization}. Here we have provided a numerical stability result for a more general case and, as a side-effect, for the particular case in Equation \eqref{eq:LHSBLaws2}. 

In this section an analysis of the discrete Lyapunov function which results from a numerical discretisation of such an analytical Lyapunov function has been discussed. An Euler scheme was applied for temporal discretisation of a split system. An upwind scheme was also applied for the spatial discretisation. The ISS-stability for such discretised systems was proved. In the next section, the results established here are applied to a linear example and the Saint-Venant model. This section endeavours to also demonstrate how values of the parameters in the Lyapunov function are delimited.
\section{Computational applications and results}\label{sec:sec03}
The results of the previous section will now be tested computationally on specific examples. The first example will be a linear hyperbolic system of balance laws with spatially-varying coefficients and boundary disturbances presented in Section \ref{sec:computeLinear}. The second example will be a Saint-Venant system of equations which will be discussed in Section \ref{sec:computeSaint-VenantEquations}. The derivation of the equilibrium and the choice of requisite parameters for such models will be discussed in detail.

\subsection{Linear hyperbolic $ 2\times 2 $ systems of balance laws with spatially-varying coefficients and boundary disturbances}\label{sec:computeLinear}	
To illustrate Theorem \ref{thm:ThmChap01-02-Part2},  we consider the following linear hyperbolic $ 2\times 2 $ system of balance laws with spatially-varying coefficients
\begin{equation}\label{eq:2by2HSBLaws}
\partial_t\begin{bmatrix} w_1\\ w_2 \end{bmatrix} + \begin{bmatrix} \lambda_1(x) & 0  \\ 0 & \lambda_2(x) \end{bmatrix} \partial_x\begin{bmatrix} w_1\\ w_2 \end{bmatrix} + \begin{bmatrix} \gamma_{11}(x) & \gamma_{12}(x) \\ \gamma_{21}(x)& \gamma_{22}(x) \end{bmatrix}\begin{bmatrix}w_1\\ w_2\end{bmatrix} = 0,\;x \in [0,l],\; t \in [0, +\infty),
\end{equation}
where the characteristic velocities, $ \lambda_{2}(x) < 0 < \lambda_{1}(x) $ are continuously differentiable, and the coefficients of the source term,  $ \gamma_{11}(x)$,  $\gamma_{12}(x)$, $\gamma_{21}(x)$, and $\gamma_{22}(x) $ are continuous on $[0, l]$ together with an initial condition 
\begin{equation}\label{eq:2by2LHSBLaws-IC}
\begin{bmatrix} w_1(x,0) \\ w_2(x,0) \end{bmatrix}  = \begin{bmatrix}  f(x)\\ g(x) \end{bmatrix},\; x \in (0,l),
\end{equation}
where $ f $ and $ g $ are smooth functions, boundary conditions with disturbances
\begin{equation}\label{eq:2by2LHSBLaws-BCs}
\begin{bmatrix} w_1(0,t) \\ w_2(l,t) \end{bmatrix}  = \begin{bmatrix}
0 & \kappa_{12} \\ \kappa_{21} & 0 \end{bmatrix}\begin{bmatrix} w_1(l,t) \\ w_2(0,t) \end{bmatrix} + \begin{bmatrix} m_1 & 0 \\ 0 & m_2 \end{bmatrix}\begin{bmatrix} b_1(t) \\ b_2(t) \end{bmatrix},\; t \in [0,+\infty),
\end{equation}
and compatibility conditions
\begin{equation}\label{eq:2by2LHSBLaws-CCs}
\begin{bmatrix} w_1(0,0) \\ w_2(l,0) \end{bmatrix}  = \begin{bmatrix}
0 & \kappa_{12} \\ \kappa_{21} & 0 \end{bmatrix}\begin{bmatrix} w_1(l,0) \\ w_2(0,0) \end{bmatrix},
\end{equation}
where  $ \kappa_{12} $, $ \kappa_{21} $, $ m_1 $ and $ m_2 $ are constant parameters. We adopt the assumption \textbf{A1-A4} for the system \eqref{eq:2by2HSBLaws} with boundary conditions \eqref{eq:2by2LHSBLaws-BCs}.

At steady-state, the system \eqref{eq:2by2HSBLaws} can be expressed as a linear system of ordinary differential equations with variable coefficients 
\begin{equation}\label{eq:2by2HSBLaws-SteadyState}
\dfrac{d}{dx}\begin{bmatrix} w_1^*(x)\\ w_2^*(x) \end{bmatrix} = \begin{bmatrix} -\frac{\gamma_{11}(x)}{\lambda_1(x)} & -\frac{\gamma_{12}(x)}{\lambda_1(x)} \\ - \frac{\gamma_{21}(x)}{\lambda_2(x)}& -\frac{\gamma_{22}(x)}{\lambda_2(x)} \end{bmatrix}\begin{bmatrix}w_1^*(x)\\ w_2^*(x)\end{bmatrix}, \; x \in [0,l],
\end{equation}
where $ w_1^*(x) $ and $ w_2^*(x) $ are non-uniform steady-state solutions. The solution of the system of ODEs \eqref{eq:2by2HSBLaws-SteadyState} may be computed by using Wronskian and Liouville's Formula or by Lagrange Method.

Based on the discussion in Section \ref{sec:sec02}, the system \eqref{eq:2by2HSBLaws} can be split and discretised together with the initial condition \eqref{eq:2by2LHSBLaws-IC}, the boundary conditions \eqref{eq:2by2LHSBLaws-BCs} and the compatibility conditions \eqref{eq:2by2LHSBLaws-CCs} as follows 
\begin{subequations}\label{eq:Disc2by2HSBLaws}
	\begin{align}
	&\begin{bmatrix} \widetilde{w_1}_{j}^{n}\\ \widetilde{w_2}_{j}^{n} \end{bmatrix} =  \begin{bmatrix} {w_1}_{j}^{n}\\{w_2}_{j}^{n} \end{bmatrix}  - {\frac{\Delta t}{\Delta x}} \begin{bmatrix}{\lambda}_{1,j-1} & 0\\0 & {\lambda_{2,j+1}} \end{bmatrix} \begin{bmatrix} {w_1}_{j}^{n}-{w_1}_{j-1}^{n}\\{w_2}_{j+1}^{n} - {w_2}_{j}^{n} \end{bmatrix}, \quad n = 0, \dots , N-1,\label{eq:Disc2by2HSBLaws01}\\
	&\begin{bmatrix} {w_1}_{j}^{n+1}\\{w_2}_{j}^{n+1} \end{bmatrix} = \begin{bmatrix} \widetilde{w_1}_{j}^{n}\\\widetilde{w_2}_{j}^{n} \end{bmatrix} - {\Delta t} \begin{bmatrix} \gamma_{11,j} & \gamma_{12,j} \\ \gamma_{21,j} & \gamma_{22,j} \end{bmatrix} \begin{bmatrix} \widetilde{w_1}_{j}^{n}\\ \widetilde{w_2}_{j}^{n} \end{bmatrix},\quad j = 0, \dots, J-1,\label{eq:Disc2by2HSBLaws02}
	\end{align}
	\begin{equation}\label{eq:Disc2by2HSBLawsIC}
	w_{1,j}^0 = f_j,\quad w_{2,j}^0 = g_j,\; j = 0, \dots,J-1,
	\end{equation}
	\begin{equation}\label{eq:Disc2by2HSBLawsBCs}
	\begin{bmatrix} {w_1}_{-1}^{n+1}\\{w_2}_{J}^{n+1} \end{bmatrix} =  \begin{bmatrix}0 & \kappa_{12}\\\kappa_{21} & 0 \end{bmatrix} \begin{bmatrix} {w_1}_{J-1}^{n+1}\\{w_2}_{0}^{n+1} \end{bmatrix} + \begin{bmatrix} m_1 & 0 \\ 0 & m_2 \end{bmatrix}\begin{bmatrix} b_1^{n+1} \\ b_2^{n+1} \end{bmatrix},\quad n = 0, \dots , N-1,
	\end{equation}
	\begin{equation}\label{eq:Disc2by2HSBLawsCCs}
	\begin{bmatrix} {w_1}_{-1}^{0}\\{w_2}_{J}^{0} \end{bmatrix} =  \begin{bmatrix}0 & \kappa_{12}\\\kappa_{21} & 0 \end{bmatrix} \begin{bmatrix} {w_1}_{J-1}^{0}\\{w_2}_{0}^{0} \end{bmatrix}.
	\end{equation}
\end{subequations}

For a fixed $ T > 0 $, we let the CFL condition hold. i.e.$\displaystyle \frac{\Delta t}{\Delta x}\max_{0 \leq j \leq J-1}\{|\lambda_{1,j}|,|\lambda_{2,j}|\} \leq 1$. By using the candidate discrete ISS-Lyapunov function \eqref{eq:DiscLyapunovfunCh01}, we analyse the discrete ISS of the discretised system \eqref{eq:Disc2by2HSBLaws}. For this reason, we give conditions for the assumptions in Theorem \ref{thm:ThmChap01-02-Part2}. These assumptions are 
\begin{itemize}
	\item[\textbf{C1:}] the matrix 
	\begin{equation*}
	{\theta}_{j} := -\begin{bmatrix} {{\lambda_1}_{j-1}}\left(\frac{{p_1}_{j+1} - {p_1}_{j}}{\Delta x} \right) + \left(\frac{{{\lambda_1}_{j}} - {{\lambda_1}_{j-1}}}{\Delta x}\right){p_1}_{j+1}   & 0 \\ 0 &  {{\lambda_2}_{j+1}}\left(\frac{{p_2}_{j} - {p_2}_{j-1}}{\Delta x} \right) +  \left(\frac{{{\lambda_2}_{j+1}} - {{\lambda_2}_{j}}}{\Delta x}\right){p_2}_{j-1} \end{bmatrix},
	\end{equation*}
	is positive definite for all $ j = 0, \dots, J-1$,
	\item[\textbf{C2:}] the matrix 
	\begin{multline*}
	M_{j} := \begin{bmatrix} {p_1}_{j} & 0 \\ 0 & {p_2}_{j} \end{bmatrix}\begin{bmatrix} {\gamma_{11}}_{j} & {\gamma_{12}}_{j} \\ {\gamma_{21}}_{j} & {\gamma_{22}}_{j} \end{bmatrix} + \begin{bmatrix} {\gamma_{11}}_{j} & {\gamma_{12}}_{j} \\ {\gamma_{21}}_{j} & {\gamma_{22}}_{j} \end{bmatrix}^{\top} \begin{bmatrix} {p_1}_{j} & 0 \\ 0 & {p_2}_{j} \end{bmatrix}\nonumber\\
	-{\Delta t} \begin{bmatrix} {\gamma_{11}}_{j} & {\gamma_{12}}_{j} \\ {\gamma_{21}}_{j} & {\gamma_{22}}_{j} \end{bmatrix}^{\top}\begin{bmatrix} {p_1}_{j} & 0 \\ 0 & {p_2}_{j} \end{bmatrix}\begin{bmatrix} {\gamma_{11}}_{j} & {\gamma_{12}}_{j} \\ {\gamma_{21}}_{j} & {\gamma_{22}}_{j} \end{bmatrix},
	\end{multline*} 
	is positive semi-definite for all $ j = 0, \dots, J-1$, and
	\item[\textbf{C3:}] the matrix 
	\begin{equation*}
	B_c:=\begin{bmatrix} {\lambda_1}_{J-1} {p_1}_{J} & 0 \\ 0 & |{\lambda_2}_{0}| {p_2}_{-1}\end{bmatrix}	- (1 + \xi)\begin{bmatrix}0 & \kappa_{12}\\\kappa_{21} & 0 \end{bmatrix}^{\top} \begin{bmatrix} {\lambda_1}_{-1} {p_1}_{0} & 0 \\ 0 & |{\lambda_2}_{J}| {p_2}_{J-1}\end{bmatrix} \begin{bmatrix}0 & \kappa_{12}\\\kappa_{21} & 0 \end{bmatrix},
	\end{equation*}
	is positive semi-definite for all $ \xi > 0 $.
\end{itemize}

The first assumption, \textbf{C1} holds true if both diagonal entries of $ \theta_{j} $ are positive for all $ j = 0, \dots, J-1$. i.e.
\begin{align*}
&{\eta_1}_{j} := \left(-\frac{{\lambda_1}_{j-1}}{{p_1}_{j}}\left(\frac{{p_1}_{j+1} - {p_1}_{j}}{\Delta x} \right) - \left(\frac{{{\lambda_1}_{j}} - {{\lambda_1}_{j-1}}}{\Delta x}\right)\frac{{p_1}_{j+1}}{{p_1}_{j}} \right) {p_1}_{j} > 0, \\
&{\eta_2}_{j}:= \left(- \frac{{\lambda_2}_{j+1}}{{p_2}_{j}}\left(\frac{{p_2}_{j} - {p_2}_{j-1}}{\Delta x} \right) - \left(\frac{{{\lambda_2}_{j+1}} - {{\lambda_2}_{j}}}{\Delta x}\right)\frac{{p_2}_{j-1}}{{p_2}_{j}} \right) {p_2}_{j} > 0,
\end{align*}
for all $ j = 0, \dots, J-1$. The second assumption, \textbf{C2} holds true if the matrix $ M_{j} $, which can be rewritten as 
\begin{equation*}
M_{j} = \begin{bmatrix} {M_{11}}_{j} & {M_{12}}_{j} \\ {M_{12}}_{j} & {M_{22}}_{j} \end{bmatrix},\;j = 0, \dots, J-1, 
\end{equation*}
where 
\begin{align*}
M_{11,j} &:= 2 {\gamma_{11}}_{j}{p_1}_{j}-{\Delta t} \left( {\gamma_{11}^{2}}_{j}{p_1}_{j} + {\gamma_{21}^{2}}_{j} {p_2}_{j}\right), \\
M_{12,j}&:=  {\gamma_{21}}_{j}{p_2}_{j}+{\gamma_{12}}_{j}{p_1}_{j}
- {\Delta t} \left( {\gamma_{11}}_{j} {\gamma_{12}}_{j}{p_1}_{j} + {\gamma_{21}}_{j}{\gamma_{22}}_{j}{p_2}_{j} \right),\\
M_{22,j}&:= 2 {\gamma_{22}}_{j}{p_2}_{j} - {\Delta t} \left( {\gamma_{12}^{2}}_{j}{p_1}_{j} + {\gamma_{22}^{2}}_{j}{p_2}_{j}\right),  
\end{align*} 
has non-negative eigenvalues, 
\begin{equation*}
\sigma_j^\pm = \frac{1}{2}\left( \left({M_{11}}_{j} + {M_{22}}_{j}\right) \pm \sqrt{\left({M_{11}}_{j} + {M_{22}}_{j}\right)^2 - 4\left({M_{11}}_{j} {M_{22}}_{j} - {M_{12}^2}_{j}\right)}\right) \geq 0,
\end{equation*}
for all  $ j = 0, \dots, J-1$. In the third assumption, \textbf{C3} the matrix $ B_c $ is rewritten as  
\begin{equation*}
B_c = \begin{bmatrix} {\lambda_1}_{J-1} {p_1}_{J} - \kappa_{21}^2 (1 + \xi)|{\lambda_2}_{J}| {p_2}_{J-1} & 0 \\ 0 & |{\lambda_2}_{0}| {p_2}_{-1} - \kappa_{12}^2 (1 + \xi){\lambda_1}_{-1} {p_1}_{0} \end{bmatrix}.
\end{equation*}
Then, the third assumption, \textbf{C3} holds if we can choose the parameters, $ \kappa_{12} $ and $ \kappa_{21} $ as 
\begin{equation*}
\kappa_{12}^2  \leq \frac{|{\lambda_2}_{0}| {p_2}_{-1}}{(1 + \xi){\lambda_1}_{-1} {p_1}_{0}}, \quad \text{and} \quad \kappa_{21}^2  \leq \frac{{\lambda_1}_{J-1} {p_1}_{J}}{(1 + \xi)|{\lambda_2}_{J}| {p_2}_{J-1}},
\end{equation*}
for all $ \xi > 0 $.

Based on the above assumptions \textbf{C1 - C3}, we conclude that the discrete $L^2-$function \eqref{eq:DiscLyapunovfunCh01} is a discrete ISS-Lyapunov function for the discretised system \eqref{eq:Disc2by2HSBLaws} and we approximate the time derivative of the ISS-Lyapunov function \eqref{eq:LyapunovfunCh01} by \eqref{eq:DtDiscLyapunovfun00} with $ \eta := \min_{0 \leq j \leq J-1} \{{\eta_1}_{j}, {\eta_2}_{j} \} $ and $ \nu  = \max \{{\lambda_1}_{-1}{p_1}_{0}m_1^2, |{\lambda_2}_{J}|{p_2}_{J-1}m_2^2 \} $. Therefore, we showed that the conditions of discrete ISS are satisfied for the steady-state $\displaystyle \begin{bmatrix} {w_1}_{j}^{n}&{w_2}_{j}^{n} \end{bmatrix}^\top \equiv 0,\; j = 0, \dots, J-1,\; n = 0, \dots, N-1 $ of the discretised system \eqref{eq:Disc2by2HSBLaws}. Moreover, the upper bound of the discrete ISS-Lyapunov function is defined by \eqref{eq:DtDiscLyapunovfun-04}.

To show the numerical analysis working, we analyse a $ 2\times 2 $ uniform linear system in the sense of Corollary \ref{col:Col02-Part2}. For this reason, we consider the system \eqref{eq:2by2HSBLaws} with uniform matrix coefficients of the form 
\begin{equation*}
\Lambda = \begin{bmatrix} 1 & 0 \\ 0 & -1 \end{bmatrix}, \quad \Gamma = \begin{bmatrix} 0.3 & -0.1 \\ -0.1 & 0.3 \end{bmatrix}, 
\end{equation*}
an initial condition of the form
\begin{equation}\label{eq:TestIC}
\begin{bmatrix} w_1(x,0) \\  w_2(x,0) \end{bmatrix} = \begin{bmatrix}
-0.5 \\ 0.5 \end{bmatrix},\; x \in (0,1),
\end{equation}
boundary conditions \eqref{eq:2by2LHSBLaws-BCs} with $ m_1 = m_2 = 1 $ and the rate of the boundary disturbance functions taken as $ b_1(t) = -b_2(t) = d(t) $, where 
\begin{equation*}
d(t) = \begin{dcases} 0.01 \sin^2(\pi t), & 0 \leq t < 5,\\
0, &  t \geq 5.
\end{dcases}
\end{equation*}
and compatibility conditions \eqref{eq:2by2LHSBLaws-CCs}. Then, the discretisation of the system with initial condition, boundary conditions and compatibility conditions are given by \eqref{eq:Disc2by2HSBLaws}. 

Let the CFL condition, $ \lambda\frac{\Delta t}{\Delta x} \leq 1 $, where $ \lambda = \max \{ \lambda_1, |\lambda_2|\} = 1$ holds for a fixed $ T > 0 $. Define an implicit discrete weight function by $ \displaystyle P_j := \text{diag}\{ {p_1}\exp(-\mu x_{j}) , {p_2}\exp(\mu x_{j}) \}, \; p_1 > 0, p_2 > 0, \mu > 0,\; j = 0, \dots, J-1 $. Thus, the assumption \textbf{C1} holds if the 
\begin{align*}
\eta =&\; \min\{\eta_1 , \eta_2 \},\\  
=&\; \min \left\lbrace \lambda_{1} \left(\frac{1 - \exp\left(-\mu {\Delta x}\right)}{\Delta x}\right), |\lambda_{2}| \left(\frac{1 - \exp\left(-\mu {\Delta x}\right)}{\Delta x}\right) \right\rbrace,\\
=&\; \alpha \left(\frac{1 - \exp\left(-\mu {\Delta x}\right)}{\Delta x}\right),\\
\geq&\;  \mu \alpha \exp\left(-\mu {\Delta x}\right),
\end{align*}
where $\alpha = \min \{\lambda_{1}, |\lambda_{2}| \} = 1$. Therefore, the decay rate of the ISS-Lyapunov function explicitly defined as $\displaystyle \eta_{\text{decayrate}} = \mu \exp\left(-\mu {\Delta x}\right) $. Beside that, if we can choose a sufficiently small $ \mu > 0 $ such that $ p_1 = p_2 = p > 0 $, then the assumption \textbf{C2} holds. In addition to the choices of parameters, we can fix $ \xi > 0 $ and choose   
\begin{align*}
|\kappa_{12}| &\leq \sqrt{\frac{|{\lambda_2}| {p_2}\exp(\mu x_{-1})}{(1 + \xi){\lambda_1} {p_1}\exp(-\mu x_{0})}} = \sqrt{\frac{1}{1 + \xi}},\; \text{and}\\
|\kappa_{21}| &\leq \sqrt{\frac{{\lambda_1} {p_1}\exp(-\mu x_{J})}{(1 + \xi)|{\lambda_2}| {p_2}\exp(\mu x_{J-1})}} = \sqrt{\frac{1}{1 + \xi}}\exp(-\mu),
\end{align*}
to show that the assumption \textbf{C3} holds. Furthermore, we define 
\begin{align*}
\nu  &= \max_{\mu} \{{\lambda_1}{p_1}\exp(-\mu x_{0})m_1^2, |{\lambda_2}|{p_2}\exp(\mu x_{J-1}) m_2^2 \}, \\
& = \max_{\mu} \{p\exp(-\mu x_{0})m_1^2, p\exp(\mu x_{J-1})\}, \\
&= p \exp(\mu x_{J-1}) = p\exp(\mu (1 - 0.5{\Delta x})). 
\end{align*}

We now take CFL = 0.75, $ T = 10 $,  $ {\Delta x} = 1/1600 $, $ {\Delta t} = 0.75/1600 $ and $ \xi = 0.125 $. Then, the decay rate is given by $ \eta = \mu \exp(-\mu {\Delta x})$. We also take $ p_1 = p_2 = 1 $ for $ \mu = 0.575 $. As a result, the control parameters are given by $|\kappa_{12}| \leq 0.9428 $ and $ |\kappa_{21}| \leq 0.5305 $. Therefore, the upper bound of the discrete ISS-Lyapunov function is defined by \eqref{eq:DtDiscLyapunovfun-04} with $ \nu = 1.7768 $. 

Hence, we compute a comparison of the discrete ISS-Lyapunov function and its upper bound for CFL = 0.75 and CFL = 1 in tables \ref{table:Convergence-01} and \ref{table:Convergence-02}, respectively.
\begin{table}[H]
	\centering
	\begin{tabular}{ccccccc}
		\hline 
		$ J $ & $ \|\mathcal{L}_{\text{up}}^n - \mathcal{L}^n \|_{L^\infty} $ & $ \|\mathcal{L}_{\text{up}}^n - \mathcal{L}^n \|_{L^2} $ &  $\mu$ & $\eta $ \\[0.5ex] 
		\hline 
		200    & 0.23286    &  0.36365   & 0.575    & 0.57335 \\[1ex]
		400    & 0.23069    & 0.36113    & 0.575    & 0.57417 \\[1ex]
		800    & 0.22918    & 0.35931    & 0.575    & 0.57459 \\[1ex]
		1600   & 0.22813    & 0.35801    & 0.575    & 0.57479 \\[0.5ex]
		\hline 
	\end{tabular} 
	\caption[]%
	{The comparison of the upper bound of Lyapunov function with discrete Lyapunov function.  Under CFL = 0.75, $ {\Delta x} = \frac1J $, $ {\Delta t} =  \frac{{\Delta x}}{\max \{\lambda_{1}, |\lambda_{2}| \}}\text{CFL}$, $  \xi = 0.125 $,  $ T = 10 $ and $ \kappa_{12} = \kappa_{21} = 0.5$.}
	\label{table:Convergence-01}
\end{table}

\begin{table}[H]
	\centering
	\begin{tabular}{ccccccc}
		\hline 
		$ J $ & $ \|\mathcal{L}_{\text{up}}^n - \mathcal{L}^n \|_{L^\infty} $ & $ \|\mathcal{L}_{\text{up}}^n - \mathcal{L}^n \|_{L^2} $ &  $\mu$ & $\eta $ \\[0.5ex] 
		\hline 
		200    & 0.23026    & 0.32884    & 0.575    & 0.57335 \\[1ex]
		400    & 0.22886    & 0.32746    & 0.575    & 0.57417 \\[1ex]
		800    & 0.2279    &  0.32645   & 0.575    & 0.57459 \\[1ex]
		1600   & 0.22723    &  0.32572   & 0.575    & 0.57479 \\[0.5ex]		\hline 
	\end{tabular} 
	\caption[]%
	{The comparison of the upper bound of Lyapunov function with discrete Lyapunov function.  Under CFL = 1, $ {\Delta x} = \frac1J $, $ {\Delta t} =  \frac{{\Delta x}}{\max \{\lambda_{1}, |\lambda_{2}| \}}\text{CFL}$, $  \xi = 0.125 $,  $ T = 10 $ and $ \kappa_{12} = \kappa_{21} = 0.5$.}
	\label{table:Convergence-02}
\end{table}  

From  Table \ref{table:Convergence-01} and \ref{table:Convergence-02} above, it can be observed that the rate of decay, $ \eta $ converges to $ \mu $ and both $ L^\infty $ and $ L^2 $ norm are steadily decaying. It must be noted that the role of the CFL is for numerical stability and it can be observed above that it does not play a significant role in the convergence of $\eta$.

\subsection{Saint-Venant equations}\label{sec:computeSaint-VenantEquations}
We consider flow of water in the presence of flow rate measurements error at the boundaries. 

One of the causes of disturbances of a flow of water along an open channel can be a measurement error at the ends of the channel. Thus, we study a flow of water along a prismatic channel with a rectangular cross-section, a length of $ l $ units and constant bottom slope. We consider boundary measurements in this flow. The model of the flow is described by Saint-Venant equations (see \cite{bastin2008using, Agu_2017}) of the form
\begin{equation}\label{eq:Saint-VenantEquations}
\begin{split}
&{\partial}_t H + {\partial}_x \left(HV\right) = 0,\\
&{\partial}_t V + {\partial}_x \left( \frac{1}{2} V^2 + gH\right) + \left(gC_f\frac{V^2}{H} - gS_b\right) = 0,\; x \in [0,l],\; t \in [0, +\infty),
\end{split}
\end{equation}
where $ H $ and $ V $ denote the depth and velocity of the water, respectively. Other constants, $ g $, $ C_f $, and $ S_b $ represent the gravitational constant, a friction parameter and the constant bottom slope of the channel, respectively. 
We set an initial condition
\begin{equation}\label{eq:SVEqnsIC}
H(x,0) = H_0(x), \quad V(x,0) =  V_0(x), \quad x \in (0,l),
\end{equation}
boundary conditions with disturbances 
\begin{equation}\label{eq:SVEqnsBCs}
V(0,t) = k_0 H(0,t) + b_1(t), \quad  V(l,t) = k_lH(l,t) + b_2(t),\quad t \in (0, \infty),
\end{equation}
and compatibility conditions with no disturbance at $ t = 0 $, 
\begin{equation}\label{eq:SVEqnsCCs}
V(0,0) = k_0 H(0,0), \quad  V(l,0) = k_lH(l,0),
\end{equation}
where $ k_0, k_l $ are boundary control parameters, and $ b_1$, $b_2$ are disturbance functions. 

We consider a sub-critical flow i.e. $ V^2 < gH $. Then, the system \eqref{eq:Saint-VenantEquations} can be written in the form of the system \eqref{eq:2by2HSBLaws} (The details of the calculation can be found in \cite{bastin2016stability}) with    
\begin{align*}
w_1&:= V - V^*(x) + (H - H^*(x))\sqrt{\frac{g}{H^*(x)}} , \\
w_2 & := V - V^*(x) - (H - H^*(x))\sqrt{\frac{g}{H^*(x)}} ,\\
\lambda_1(x)&:= V^*(x) + \sqrt{gH^*(x)}, \quad  \lambda_2(x):= V^*(x) - \sqrt{gH^*(x)},\\
\gamma_{11}(x)&:= \frac{3}{4}\frac{g}{H^*(x)}\left(\frac{S_bH^*(x) - C_f{V^*}^2(x)}{\lambda_1(x)} \right) + \frac{gC_f{V^*}^2(x)}{2H^*(x)}\left(\frac{2}{V^*(x)} - \frac{1}{\sqrt{gH^*(x)}}\right),\\
\gamma_{12}(x)&:= \frac{1}{4}\frac{g}{H^*(x)}\left(\frac{S_bH^*(x) - C_f{V^*}^2(x)}{\lambda_1(x)} \right) + \frac{gC_f{V^*}^2(x)}{2H^*(x)}\left(\frac{2}{V^*(x)} + \frac{1}{\sqrt{gH^*(x)}}\right),\\
\gamma_{21}(x)&:= \frac{1}{4}\frac{g}{H^*(x)}\left(\frac{S_bH^*(x) - C_f{V^*}^2(x)}{\lambda_2(x)} \right) + \frac{gC_f{V^*}^2(x)}{2H^*(x)}\left(\frac{2}{V^*(x)} - \frac{1}{\sqrt{gH^*(x)}}\right),\\
\gamma_{22}(x)&:= \frac{3}{4}\frac{g}{H^*(x)}\left(\frac{S_bH^*(x) - C_f{V^*}^2(x)}{\lambda_2(x)} \right) + \frac{gC_f{V^*}^2(x)}{2H^*(x)}\left(\frac{2}{V^*(x)} + \frac{1}{\sqrt{gH^*(x)}}\right),
\end{align*}
where $ H^*(x) $, $ V^*(x) $ is an equilibrium solution.
Also, the initial condition \eqref{eq:SVEqnsIC}, the boundary conditions \eqref{eq:SVEqnsBCs}, and the compatibility conditions \eqref{eq:SVEqnsCCs} are expressed as \eqref{eq:2by2LHSBLaws-IC}, \eqref{eq:2by2LHSBLaws-BCs}, and \eqref{eq:2by2LHSBLaws-CCs}, respectively with 
$ f(x):= w_1(x,0)$, $ g(x):= w_2(x,0) $, $ \kappa_{12}:= \frac{k_0\sqrt{\frac{H^*(0)}{g}} - 1}{1 + k_0 \sqrt{\frac{H^*(0)}{g}}} \neq 1$, $ \kappa_{21}:= \frac{k_l\sqrt{\frac{H^*(l)}{g}} - 1}{1 + k_l \sqrt{\frac{H^*(l)}{g}}} \neq 1$, $ m_1:= 1 - \kappa_{12} $ and $ m_2:= 1 - \kappa_{21}$. 

For a numerical analysis and computations, we take an example from \cite{diagne2017backstepping}. Thus, a constant steady-state solution, $ H^*(x) = 2, V^*(x) = 3,\; x \in [0,1] $ is considered. The parameters are given by $ g = 9.81$, $C_f = 0.1 $ and $ S_b = 0.0459 $, and initial condition defined by $ H(x,0) = 2.5 $, $ V(x,0) = 4\sin(\pi x) $ for  $ x \in [0,1] $.  

Therefore, $ \lambda_1 = 7.4294$, $\lambda_2 = -1.4294$, $\gamma_{11}(x) = \gamma_{21}(x) = 0.0992$ and $\gamma_{12}(x) = \gamma_{22}(x) = 0.2008$ for all $ x \in [0,1] $. We set an initial condition $ w_1(x,0) := -1.8926 + 4\sin(\pi x) $, $w_2(x,0) = -4.1074 + 4\sin(\pi x), \; x \in (0,1) $. The rate of the boundary disturbance functions taken as $ b_1(t) = -b_2(t) = d(t) $, where 
\begin{equation*}
d(t) = \begin{dcases} 0.01 \sin^2(\pi t), & 0 \leq t < 5,\\
0, &  t \geq 5.
\end{dcases}
\end{equation*}

We now take CFL = 0.75, $ T = 10 $,  $ {\Delta x} = 1/1600 $, $ {\Delta t} = 0.75/1600 \lambda $, where  $ \lambda = \max \{ \lambda_1, |\lambda_2|\} = 7.4294$. Define a discrete weight function $ P_j = \text{diag} \{p_1 e^{-\mu x_{j}}, p_2 e^{\mu x_{j}} \}$, $j = 0, \dots, J-1. $ Then, the decay rate is given by $ \eta = \mu \alpha \exp(-\mu {\Delta x})$, where $ \alpha = \min \{ \lambda_1, |\lambda_2|\} = 1.4294$. A sufficiently small value of $\mu$ can be chosen such that $p_1\gamma_{12} = p_2\gamma_{21}$. Thus, $p_1 = \gamma_{21} = 0.0992$ and $p_2 = \gamma_{12} = 0.2008$. We fix $ \xi = 0.125 $, then the control parameters are given by $|\kappa_{12}| < 0.5884$ and $|\kappa_{21}| < 1.5108\exp(-\mu)$. With the choice of boundary control parameters $ \kappa_{12} $, $ \kappa_{21} $, the coefficients of boundary disturbance functions can be obtained as $ m_1 = 1 - \kappa_{12} $, $ m_2 = 1 -\kappa_{21} $. Therefore, the upper bound of the discrete ISS-Lyapunov function is defined by \eqref{eq:DtDiscLyapunovfun-04} with 
\begin{equation*}
\nu  = \max_{\mu} \{{\lambda_1}{p_1}\exp(-\mu x_{0})m_1^2, |{\lambda_2}|{p_2}\exp(\mu x_{J-1}) m_2^2 \}. 
\end{equation*}

\begin{figure}[h!]
	\centering
	\includegraphics{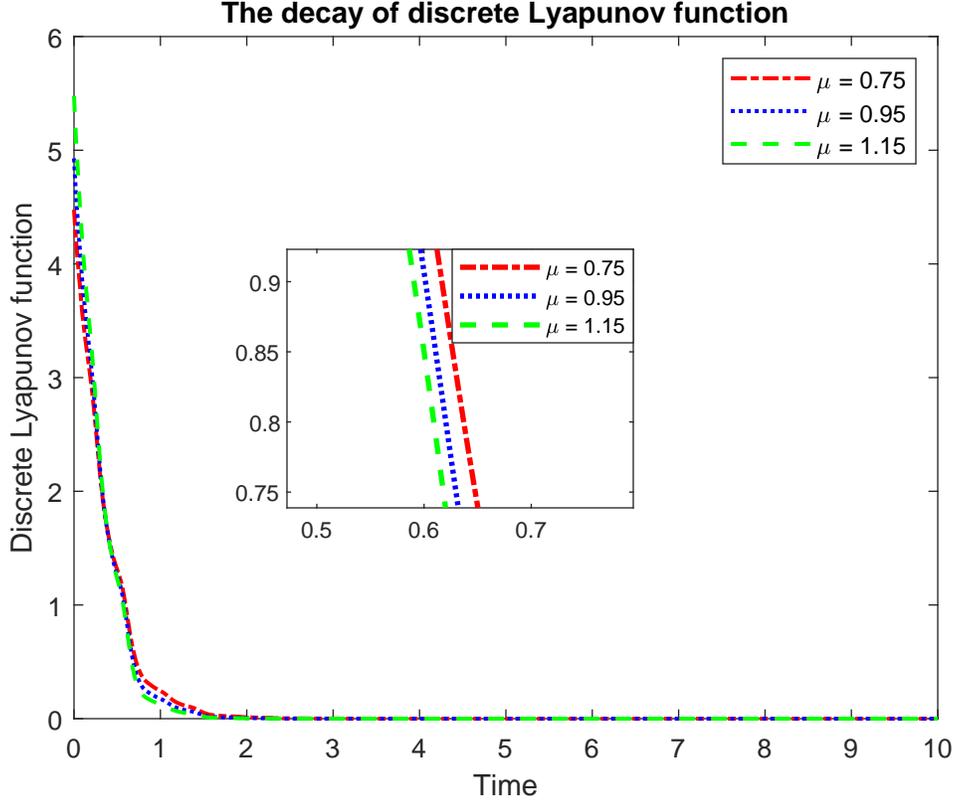}
	\caption{The decay of Lyapunov function for Saint-Venant equations. The choice of parameters are $ p_1 =0.0992,  p_2 = 0.2008 $, $ \kappa_{12} = 0.5 $, $ \kappa_{21} = 1.5\exp(-\mu)$, $ m_1 = 1 - \kappa_{12} $ and $ m_2 = 1 -\kappa_{21} $  with  $ l = 1 $, $ J = 1600 $ and $ T = 10$ under CFL = 0.75.}\label{fig:SVEqnswithBD}
\end{figure}

In Figure \ref{fig:SVEqnswithBD}, it can be observed that the three nearly indistinguishable curves which are obtained for different values of $ \mu > 0 $ converge to 0 asymptotically in time. This shows the decay of the ISS-Lyapunov function in the presence of boundary disturbance. Hence, in the sense of the definition of discrete ISS, the steady-state $ W_j^n \equiv 0,\; j = 0, \dots, J-1,\; n = 0, \dots, N-1 $ of the discretised system with the discretised boundary conditions is discrete ISS in $L^2-$norm with respect to discrete disturbance function $ b^n,\;n=0,\dots, N-1 $. 

Similar computations were also applied to the isothermal Euler equations for which Condition \textbf{C2} does not hold. We have taken an example in  \cite{gersterdiscretized}, $ \rho^*(0) = 3, q^*(x) = 0.2,\; x \in [0,1] $ with the parameters given by $ a =1, \frac{f}{D} = 1 $. Thus \[{\rho}^*(x) = \frac{3}{\exp \left(\frac{\text{LambertW}(-1, -225\exp(x-225))}{2}-\frac{x}{2} + \frac{225}{2}\right)}.\] 

We considered the system \eqref{eq:2by2HSBLaws}, the initial condition \eqref{eq:2by2LHSBLaws-IC}, the boundary conditions \eqref{eq:2by2LHSBLaws-BCs} and the compatibility conditions \eqref{eq:2by2LHSBLaws-CCs} with 
\begin{align*}
\frac{{q^*}(x)}{\rho^*(x)} - a =:& \lambda_2(x) < 0 < \lambda_1(x) := \frac{{q^*}(x)}{\rho^*(x)} + a, \\
\gamma_{11}(x) =&\; -\frac{1}{2a} \left( \lambda_2(x)\dfrac{d}{dx} \lambda_1(x) + \lambda_1(x) \dfrac{d}{dx} \lambda_2(x) + \frac{f}{D} \frac{{q^*}^2(x)}{2{\rho^*}^2(x)} \right)\\
&\;- \frac{\lambda_1(x)}{2a} \left(\frac{2q^*(x)}{{\rho^*}^2(x)} - \frac{f}{D} \frac{q^*(x)}{{\rho^*}(x)}\right) + \frac{1}{2a} \dfrac{d}{dx} \lambda_2(x),\\
\gamma_{12}(x) =&\; \frac{1}{2a} \left( \lambda_2(x)\dfrac{d}{dx} \lambda_1(x) + \lambda_1(x) \dfrac{d}{dx} \lambda_2(x) + \frac{f}{D} \frac{{q^*}^2(x)}{2{\rho^*}^2(x)} \right)\\
&\;+ \frac{\lambda_2(x)}{2a} \left(\frac{2q^*(x)}{{\rho^*}^2(x)} - \frac{f}{D} \frac{q^*(x)}{{\rho^*}(x)}\right) - \frac{1}{2a} \dfrac{d}{dx} \lambda_2(x),\\
\gamma_{21}(x) =&\; -\frac{1}{2a} \left( \lambda_2(x)\dfrac{d}{dx} \lambda_1(x) + \lambda_1(x) \dfrac{d}{dx} \lambda_2(x) + \frac{f}{D} \frac{{q^*}^2(x)}{2{\rho^*}^2(x)} \right)\\
&\;- \frac{\lambda_1(x)}{2a} \left(\frac{2q^*(x)}{{\rho^*}^2(x)} - \frac{f}{D} \frac{q^*(x)}{{\rho^*}(x)}\right) + \frac{\lambda_1(x)}{2a} \dfrac{d}{dx} \lambda_1(x),\\
\gamma_{22}(x) =&\; \frac{1}{2a} \left( \lambda_2(x)\dfrac{d}{dx} \lambda_1(x) + \lambda_1(x) \dfrac{d}{dx} \lambda_2(x) + \frac{f}{D} \frac{{q^*}^2(x)}{2{\rho^*}^2(x)} \right)\\
&\;+ \frac{\lambda_2(x)}{2a} \left(\frac{2q^*(x)}{{\rho^*}^2(x)} - \frac{f}{D} \frac{q^*(x)}{{\rho^*}(x)}\right) - \frac{\lambda_2(x)}{2a} \dfrac{d}{dx} \lambda_1(x),
\end{align*}
$ f(x) = g(x) = \cos(2\pi x), \; x \in (0,1) $ and the rate of the boundary disturbance functions taken as $ b_1(t) = -b_2(t) = d(t) $, where 
\begin{equation*}
d(t) = \begin{dcases} 0.01 \sin^2(\pi t), & 0 \leq t < 5,\\
0, &  t \geq 5.
\end{dcases}
\end{equation*}  
Since $ \gamma_{11}(x) > 0 $, $ \gamma_{12}(x) > 0 $, $ \gamma_{21}(x) > 0 $ but $ \gamma_{22}(x) < 0 $ for all $ x \in [0, 1] $, the matrix $ M_j $ in Condition \textbf{C2} cannot be positive semi-definite. Therefore, it cannot be guaranteed that the discrete $ L^2- $function defined by \eqref{eq:DtDiscLyapunovfun-04} is the discrete ISS-Lyapunov function. 

\section{Conclusion}\label{sec:sec04}

In this paper, we presented the discretisation of a linear hyperbolic system of balance laws with boundary disturbance. For numerical discretisation, we used a finite volume method. Specifically, we used upwind scheme and time splitting method. We also discretised an $ L^2-$ISS-Lyapunov function to investigate conditions for ISS of the discretised system. Finally, the result was applied to a linear problem and a relevant physical problem: Saint-Venant equations and numerical simulations are computed in order to test the results and compare with analytical results. We also established that for the isothermal Euler equations, one of the conditions required for ISS are not satisfied hence the result in this paper may not hold. The properties that have been proved analytically can also be established computationally.

This work leaves more questions open. There is need to analyse Lyapunov functions for nonlinear differential equations. Analysis of numerical artefacts such as numerical viscosity need to be carefully examined. Such numerical artefacts may have an influence on the rate of convergence of the discrete results.




\end{document}